\documentclass{article}
\usepackage[round]{natbib}
\usepackage{amssymb}
\usepackage{amsmath}
\usepackage{amsthm}
\usepackage{graphicx}
\usepackage{verbatim}
\usepackage{framed}
\usepackage{tikz,subfigure}
\usepackage{multicol}
\usepackage[colorlinks=true,citecolor=black,linkcolor=black,urlcolor=blue]{hyperref}
\usepackage[margin=4cm]{geometry}
\usepackage{linguex}
\usepackage{relsize}
\usepackage{proof}

\tikzstyle{index gray}=[inner sep=3pt, black, circle, fill=lightgray]
\tikzstyle{opaque}=[fill=gray,fill opacity=.1]

\DeclareFontFamily{U}{MnSymbolC}{}
\DeclareSymbolFont{MnSyC}{U}{MnSymbolC}{m}{n}
\DeclareMathSymbol{\diamondplus}{\mathbin}{MnSyC}{"7C}
\DeclareFontShape{U}{MnSymbolC}{m}{n}{
    <-6>  MnSymbolC5
   <6-7>  MnSymbolC6
   <7-8>  MnSymbolC7
   <8-9>  MnSymbolC8
   <9-10> MnSymbolC9
  <10-12> MnSymbolC10
  <12->   MnSymbolC12}{}

\renewcommand{\phi}{\varphi}
\def\disji{\rotatebox[origin=c]{-90}{$\!{\geqslant}$}}
\newcommand{\lori}{\,\disji\,}
\def\DISJI{\ensuremath{\mathlarger{\mathbin{{\setminus}\mspace{-5mu}{\setminus}}\hspace{-0.33ex}/}}\xspace}

\newcommand{\Lori}{\,\DISJI\,}

\newcommand{\ibox}{\boxplus}
\newcommand{\idia}{\diamondplus}
\newcommand{\kbox}{\square}
\newcommand{\R}{\mathcal{R}}
\newcommand{\vdasheq}{\dashv\vdash}
\newcommand{\Togamma}{\Rrightarrow_\Gamma}
\newcommand{\Left}{\textsf{L}}
\newcommand{\Right}{\textsf{R}}
\newcommand{\counterfactual}{\ensuremath{%
  \Box\kern-1.5pt
  \raise1pt\hbox{$\mathord{\rightarrow}$}}}
  
 \newcommand{\inl}{\textsf{INL}}

\renewcommand{\aa}{\alpha}
\newcommand{\bb}{\beta}
\renewcommand{\gg}{\gamma}

\newcommand{\C}{\ensuremath{\mathcal{C}}}
\newcommand{\inqb}{\textsf{InqB}}
\newcommand{\inqml}{\textsf{InqML}}
\newcommand{\inqnl}{\ensuremath{\textsf{InqML}_{\To}}}
\newcommand{\md}{\text{md}}
\newcommand{\To}{\Rrightarrow}
\newcommand{\modeq}{\equiv}

\renewcommand{\L}{\ensuremath{\mathcal{L}}}

\renewcommand{\P}{\ensuremath{\mathcal{P}}}

\renewcommand\>{\rangle}
\newcommand{\MM}{\ensuremath{\mathcal{M}}}
\newcommand{\M}{\ensuremath{\mathsf{M}}}

\theoremstyle{definition}
\newtheorem{theorem}{Theorem}[section]
\newtheorem{definition}[theorem]{Definition}
\newtheorem{proposition}[theorem]{Proposition}
\newtheorem{observation}[theorem]{Observation}

\newtheorem{corollary}[theorem]{Corollary}

\newtheorem{lemma}[theorem]{Lemma}

    \newtheorem{remark}[theorem]{Remark}
    \newtheorem{example}[theorem]{Example}

\title{Inquisitive Neighborhood Logic\thanks{The present paper subsumes and significantly extends a paper published in the proceedings of Advances in Modal Logic \citep{Ciardelli:22aiml}. 
Funding from the German Science Foundation (DFG, Project number 446711878) and the European Research Council (ERC, Grant agreement number 101116774) is gratefully acknowledged.}}
\author{Ivano Ciardelli}
\date{\today}

\begin{document}

\maketitle

\begin{abstract}

\noindent
We explore an inquisitive modal logic designed to reason about neighborhood models. This logic is based on an inquisitive strict conditional operator $\To$, which quantifies over neighborhoods, and which can be applied to both statements and questions. In terms of this operator we also define two unary modalities $\ibox$ and $\idia$, which act respectively as a universal and existential quantifier over neighborhoods.
We prove that the expressive power of this logic matches the natural notion of bisimilarity in neighborhood models. We show that certain fragments of the language are invariant under certain modifications of the set of neighborhoods, and use this to show that $\To$ is not definable from $\ibox$ and $\idia$, and that questions embedded under $\To$ are indispensable. We provide a sound and complete axiomatization of our logic, both in general and in restriction to some salient frame classes, and discuss the relations between our logic and other modal logics interpreted over neighborhood models.

\end{abstract}

\section{Introduction}

A neighborhood model is a modal structure where each world $w$ is associated with a set $\Sigma(w)$ of neighborhoods, where each neighborhood is a set of worlds. As an illustration, Figure~\ref{fig:1} shows the neighborhoods of three worlds, $w_1,w_2,w_3$: the neighborhoods of $w_1$ are the four singletons, $\{w_{pq}\},\{w_{p\overline q}\},\{w_{\overline pq}\},\{w_{\overline{pq}}\}$, where the subscript on a world represents the corresponding propositional valuation; the neighborhoods of $w_2$ are the same four singletons as well as their union, $\{w_{pq},w_{p\overline q},w_{\overline pq},w_{\overline{pq}}\}$; finally, the neighborhoods of  $w_3$ include, in addition, also the sets $\{w_{pq},w_{p\overline q}\}$ and $\{w_{\overline pq},w_{\overline{pq}}\}$.

\begin{figure}[t]
\centering
\label{fig:inquisitive states}
\subfigure[$\Sigma(w_1)$]{\label{fig:a}
 \begin{tikzpicture}[>=latex,scale=0.55]

 \draw[opaque,rounded corners] (-1.7,1.7) rectangle (-.3,.3);
 \draw[opaque,rounded corners] (1.7,1.7) rectangle (.3,.3);
 \draw[opaque,rounded corners,yshift=-2cm] (-1.7,1.7) rectangle (-.3,.3);
 \draw[opaque,rounded corners,yshift=-2cm] (1.7,1.7) rectangle (.3,.3);

 \draw (2,2) rectangle (2,2); 
 \draw (-2,-2) rectangle (-2,-2); 

 \draw (-1,1) node[index gray] (yy) {\footnotesize $pq$};
 \draw (1,1) node[index gray] (yn) {\footnotesize$p\overline q$};
 \draw (-1,-1) node[index gray] (yn) {\footnotesize$\overline pq$};
 \draw (1,-1) node[index gray] (yn) {\footnotesize$\overline{pq}$};

	\end{tikzpicture}
}
\hspace{.5in}
\subfigure[$\Sigma(w_2)$]{\label{fig:b}
 \begin{tikzpicture}[>=latex,scale=0.55]

 \draw[opaque,rounded corners] (-1.7,1.7) rectangle (-.3,.3);
 \draw[opaque,rounded corners] (1.7,1.7) rectangle (.3,.3);
 \draw[opaque,rounded corners,yshift=-2cm] (-1.7,1.7) rectangle (-.3,.3);
 \draw[opaque,rounded corners,yshift=-2cm] (1.7,1.7) rectangle (.3,.3);

 \draw[opaque,rounded corners] (-2,2) rectangle (2,-2);

 \draw (2,2) rectangle (2,2); 
 \draw (-2,-2) rectangle (-2,-2); 

 \draw (-1,1) node[index gray] (yy) {\footnotesize $pq$};
 \draw (1,1) node[index gray] (yn) {\footnotesize$p\overline q$};
 \draw (-1,-1) node[index gray] (yn) {\footnotesize$\overline pq$};
 \draw (1,-1) node[index gray] (yn) {\footnotesize$\overline{pq}$};

	\end{tikzpicture}
}
\hspace{.5in}
\subfigure[$\Sigma(w_3)$]{\label{fig:c}
 \begin{tikzpicture}[>=latex,scale=0.55]

 \draw[opaque,rounded corners] (-1.7,1.7) rectangle (-.3,.3);
 \draw[opaque,rounded corners] (1.7,1.7) rectangle (.3,.3);
 \draw[opaque,rounded corners,yshift=-2cm] (-1.7,1.7) rectangle (-.3,.3);
 \draw[opaque,rounded corners,yshift=-2cm] (1.7,1.7) rectangle (.3,.3);

 \draw[opaque,rounded corners] (-1.85,1.85) rectangle (1.85, .15);
 \draw[opaque,rounded corners,yshift=-2cm] (-1.85,1.85) rectangle (1.85, .15);

 \draw[opaque,rounded corners] (-2,2) rectangle (2,-2);

 \draw (1.9,1.9) rectangle (1.9,1.9); 
 \draw (-1.9,-.9) rectangle (-1.9,-.9); 

 \draw (-1,1) node[index gray] (yy) {\footnotesize $pq$};
 \draw (1,1) node[index gray] (yn) {\footnotesize$p\overline q$};
 \draw (-1,-1) node[index gray] (yn) {\footnotesize$\overline pq$};
 \draw (1,-1) node[index gray] (yn) {\footnotesize$\overline{pq}$};

	\end{tikzpicture}
}
\caption{Sets of neighborhoods associated with three worlds.}
\label{fig:1}
\end{figure}

Like Kripke models, neighborhood models can be given many concrete interpretations. Among others, they have been used to interpret logics of evidence and evidence-based belief \citep{Benthem:14}, and to define an inquisitive generalization of epistemic logic in which agents are \emph{curious}, i.e., they do not just possess certain information, but also entertain certain questions \citep{CiardelliRoelofsen:15idel}. A special class of neighborhood models, \emph{sphere models}, provides the most popular semantics for counterfactual conditionals \citep{Lewis:73} and is standardly employed to model rational belief revision \citep{Board:04,BaltagSmets:06}. Another special class, \emph{topological models}, provides a natural semantics for spatial reasoning \citep{Aiello:03}. 
Furthermore, neighborhood models emerge naturally in the context of reasoning about the abilities of individuals and groups in a strategic setting \citep{Brown:88,Pauly:02}. 

In light of these many interpretations, it is not surprising that there is a long tradition of work using neighborhood models to provide semantics for modal logic \citep[see][for an overview]{Pacuit:17}. In particular, two different neighborhood-based semantics for standard modal logic are well-established in the literature. According to the so-called Scott-Montague semantics \citep{Scott:70,Montague:70universal}, $\Box\aa$ is true at a world $w$ if the truth-set $|\aa|_M$, the set of worlds in the model $M$ where $\aa$ is true, is a neighborhood of $w$.
\begin{description}
\item[Scott-Montague semantics:] $M,w\models\Box\aa\iff |\aa|_M\in\Sigma(w)$.
\end{description}

\noindent
According to a different approach, rooted in the topological semantics of \cite{Tarski:38} and \cite{McKinseyTarski:44}, 
 \citep[see][for modern formulations]{Chellas:80,Brown:88,Aiello:03,Hansen:03}, $\Box\aa$ is true at $w$ if $\aa$ is true everywhere in some neighborhood of $w$. We refer to this as the $\exists\forall$-semantics. 
 
\begin{description}
\item[$\exists\forall$ semantics:] $M,w\models\Box\aa\iff \exists s\in\Sigma(w)\text{ s.t. }\forall v\in s: M,v\models\aa$.
\end{description}

\noindent
The two clauses coincide in the context of \emph{upward-monotone} neighborhood models, i.e., models in which every superset of a neighborhood is itself a neighborhood, but they come apart in the general case.

As it happens, however, the standard modal language with either of the above clauses 
does not allow us to express much about the configurations that arise in the neighborhoods of a given world. For instance, the situation at worlds $w_1,w_2$, and $w_3$ above is clearly different in some salient respects: for instance, every neighborhood of $w_1$ settles whether $p$ is true or false (i.e., the truth value of $p$ is constant within each neighborhood) while this is not the case for $w_2$ and $w_3$; moreover, every neighborhood of $w_2$ that settles whether $p$ also settles whether $q$, whereas this is not the case for $w_3$. Surprisingly, these simple facts cannot be expressed in the standard modal language with either of the above semantics.

A way to bring the problem more sharply into focus is to think about the notion of \emph{behavioral equivalence} in neighborhood models.  %
Intuitively, it seems natural to regard two worlds as behaviorally equivalent if they validate the same atomic formulas and  every neighborhood of the one is behaviorally equivalent to some neighborhood of the other; in turn, it seems natural to regard two neighborhoods as behaviorally equivalent when each world in the one is behaviorally equivalent to some world in the other. This leads naturally to a notion of \emph{bisimilarity} for neighborhood models, analogous to the standard notion of bisimilarity that plays a key role in the setting of relational semantics. 

Given the usefulness of standard modal logic in describing bisimulation-invariant properties of relational models, it seems worthwhile to develop a modal logic to describe the bisimulation-invariant properties of neighborhood models. The expressive power of this logic should match the notion of neighborhood bisimilarity in the same sense in which the expressive power of standard modal logic matches the standard notion of bisimilarity. In particular, bisimilarity should guarantee modal equivalence, and the opposite should hold as well in the context of finite models over a finite set of atoms.

Neither of the above interpretations of the standard modal language satisfies this desideratum. The Scott-Montague semantics is not bisimulation-invariant, since it is inherently non-local: the truth of $\Box\aa$ at a world $w$ depends on the truth value of the argument $\aa$ across the entire model, and not only locally in the neighborhoods of $w$; as a consequence, the truth of $\Box\aa$ may fail to be preserved even under the simple operation of taking the disjoint union of two models. The $\exists\forall$-semantics, by contrast, is bisimulation-invariant, but it fails to distinguish non-bisimilar worlds even in a finite setting. For instance, the three worlds in Fig.~\ref{fig:1} are not bisimilar to one another, and the difference shows up already after a single modal transition: for instance, $w_3$ has a neighborhood, $\{w_{pq},w_{p\overline q}\}$, which is not behaviorally equivalent to any neighborhood of $w_1$ or $w_2$; nevertheless, it is easy to see (cf.\ also Proposition \ref{prop:monotonicity} below) that these three worlds satisfy exactly the same modal formulas under the $\exists\forall$-semantics.

Recently, a modal logic that satisfies the above desiderata has been proposed by \cite{Benthem:17}: their approach, called \emph{instantial neighborhood logic}, is based on modal formulas of the form $\Box(\bb_1,\dots,\bb_n;\aa)$, where $n\ge 0$, which are true at a world $w$ if there is a neighborhood of $w$ where $\aa$ is true everywhere, and each $\beta_i$ is true somewhere. 

In this paper, we explore a very different way to set up a modal logic with the desired level of expressive power. Our approach is based on the framework of inquisitive logic \citep{Ciardelli:18book,Ciardelli:23book}, which means that our logical language contains not only formulas regimenting statements, as usual, but also formulas regimenting questions. For instance, our language includes, next to an atom $p$ formalizing the statement ``that $p$'', also the formula $?p$, formalizing the question ``whether $p$''. To the standard language of inquisitive propositional logic, our logic adds a binary modality $\To$, read `yields'. This modality is essentially a strict conditional operator, except for the fact that it quantifies not over worlds, but over neighborhoods; this is possible since in inquisitive logic formulas can be evaluated directly relative to a set of worlds, in terms of a relation called \emph{support}.\footnote{Inquisitive modal logics share this feature with other recently developed modal logics whose semantics is based on sets of worlds, called ``information states'' or ``teams''.  These include modal dependence and independence logics \citep{Vaananen:08,Hella:13wollic,Kontinen:17}, the bilateral state-based modal logic of \cite{Aloni:22}, and the logic of epistemic modals of \cite{HawkeSteinert-Threlkeld:21}.}
 Using the same symbol $\models$ both for truth at a world (on the left) and for support at a set of worlds (on the right), the semantics of our modality can be stated as follows:

\begin{description}
\item[Semantics of `yields':] $M,w\models\phi\To\psi\;\iff\; \forall s\in\Sigma(w): M,s\models\phi\text{ implies }{M,s\models\psi}$.
\end{description}
Crucially, the modality $\To$ may be applied not only to statements, but also to questions. For instance, the claim ``every neighborhood settles whether $p$'' (true at $w_1$ but not at $w_2,w_3$ in Fig.~\ref{fig:1}) can be expressed by the formula ${\top\To{?p}}$, whereas the claim ``every neighborhood that settles whether $p$ also settles whether $q$'' (true at $w_1,w_2$ but not at $w_3$) can be expressed as ${{?p}\To{?q}}$. Thus, by means of modal formulas that involve embedded questions, the worlds $w_1,w_2,$ and $w_3$ from Fig.~\ref{fig:1} can be distinguished. 

In terms of the modality $\To$, we may define two unary modalities, $\ibox$ (`window') and $\idia$ (`kite'), which function respectively as universal and existential quantifier over neighborhoods. 
Previous work on inquisitive modal logic \citep[see, a.o.,][]{CiardelliRoelofsen:15idel,Ciardelli:16,CiardelliOtto:18,MeissnerOtto:22,PuncocharSedlar:21epistemic,MaricPerkov:24}  has focused on a language based on the modality $\ibox$.  
As \cite{CiardelliOtto:18} have shown, this language is expressively adequate for bisimulation-invariant properties in the setting of \emph{downward-monotone} neighborhood models, i.e., models where every subset of a neighborhood is itself a neighborhood. The focus on downward-monotone models in this literature was motivated by the particular interpretation of neighborhood models that first drove the development of inquisitive modal logic, namely, the \emph{inquisitive-epistemic} interpretation \citep{CiardelliRoelofsen:15idel}. However, as discussed above, neighborhood models admit many interpretations; if one wants to view inquisitive modal logic as a tool to reason about neighborhood models in general, a language based on $\ibox$ is no longer expressively adequate; for instance, this language would not be able distinguish worlds $w_2$ and $w_3$ in Fig.~\ref{fig:1} (in fact, as we shall see, even a language with both $\ibox$ and $\idia$ cannot distinguish these worlds, cf.\ Proposition \ref{prop:monotonicity}). The logic explored in this paper can thus be seen as a generalization of standard inquisitive modal logic which is suitable for reasoning about neighborhood models which are not necessarily downward-monotone. This makes it possible to explore a much broader spectrum of concrete interpretations for inquisitive modal logic. 

In this paper, we explore in detail two aspects of the inquisitive modal logic $\inqml_\To$: its expressive power, and its axiomatization.

With respect to expressive power, we will show that our logic matches the natural notion of bisimulation for neighborhood models in exactly the same sense in which standard modal logic matches the notion of bisimulation relative to Kripke models. In a slogan:
$$\inqnl:\text{neighborhood bisimilarity}=\textsf{ML}:\text{Kripke bisimilarity}$$
In particular, over a finite set of atoms, a property of worlds is definable in \inqnl\ just in case it is invariant under some finite degree of bisimilarity. Moreover, two worlds in a finite model are bisimilar if and only if they agree on all modal formulas.

With respect to axiomatization, we will show that the logic of $\To$ is characterized by four simple logical principles: (i) if $\phi$ entails $\psi$, then $\phi\To\psi$ is a logical validity; (ii) $\To$ is transitive; (iii) if an antecedent $\phi$ yields two consequents, it yields their conjunction; (iv) if a consequent $\psi$ is yielded by two antecedents, it is yielded by their inquisitive disjunction. Moreover, we prove that certain salient conditions on neighborhood frames are captured by natural modal axioms (though, in other cases, the problem of finding a suitable axiom remains open).

The paper is structured as follows: Section~\ref{sec:background} covers some preliminaries on inquisitive propositional logic and the general architecture of inquisitive modal logics; Section~\ref{sec:inqnl} introduces the modal logic \inqnl\ and illustrates it with some examples; Section~\ref{sec:fragments} shows that certain fragments of \inqnl\ are invariant under certain modifications of the set of neighborhoods; this is used to prove that $\To$ is not definable from $\ibox$ and $\idia$ and that questions embedded on the right of $\To$ are indispensable; Section \ref{sec:expressive power} introduces the notion of bisimilarity for neighborhood models and characterizes the expressive power of \inqnl\ in terms of this notion; Section \ref{sec:axiomatization} provides an axiomatization of \inqnl\ and establishes its strong completeness; Section \ref{sec:frame conditions} extends this result by giving canonical axioms for some salient frame conditions; Section \ref{sec:empty} discusses the treatment of the empty neighborhood; Section \ref{sec:inl} relates \inqnl\ to the instantial neighborhood logic of \cite{Benthem:17}; Section~\ref{sec:conclusion} outlines some directions for future work.

\section{Background}
\label{sec:background}

In this section, we briefly cover some relevant background on inquisitive propositional logic, and on the general features of inquisitive modal logics. For a more thorough introduction to the topic, and for proofs of the results mentioned in this section, the reader is referred to \cite{Ciardelli:23book}.

\subsection{Propositional inquisitive logic}

The standard system of propositional inquisitive logic, \inqb, is a conservative extension of classical propositional logic with a question-forming disjunction, $\lori$, called \emph{inquisitive disjunction}. The language $\L$ of \inqb\ is given by the following definition, where $p$ is an atomic sentence (or `atom') drawn from a fixed set $\P$:
$$\phi\;:=\;p\mid\bot\mid(\phi\land\phi)\mid(\phi\to\phi)\mid(\phi\lori\phi)$$
Other connectives are defined as syntactic shorthands, as follows:
\begin{multicols}{2}
\begin{itemize}
\item $\neg\phi:=\phi\to\bot$
\item $\phi\lor\psi:=\neg(\neg\phi\land\neg\psi)$
\item $\phi\leftrightarrow\psi:=(\phi\to\psi)\land(\psi\to\phi)$
\item $?\phi:=\phi\lori\neg\phi$
\end{itemize}
\end{multicols}
\noindent
The $\lori$-free fragment of \L\ can be identified with the standard language of propositional logic. To mark the contrast with inquisitive disjunction, we refer to $\lor$ as \emph{declarative disjunction}. Intuitively, whereas the declarative disjunction $p\lor\neg p$ regiments the disjunctive statement \emph{that $p$ is true or false} (a tautology) the inquisitive disjunction $p\lori\neg p$ regiments the yes-no question \emph{whether $p$ is true or false} (whence the abbreviation $?p$).

Models for \inqb\ are simply pairs $M=\<W,V\>$ consisting of a non-empty universe $W$ of possible worlds and a valuation function $V:W\times\P\to\{0,1\}$ which assigns a truth-value $V(w,p)$ to an atom $p$ relative to a world $w$. 
As mentioned in the introduction, however, the semantics for \inqb\ is not given by a definition of truth at a possible world, but rather by a definition of \emph{support} relative to an \emph{information state}, defined as a set of worlds $s\subseteq W$. The notion of support is defined by the following inductive clauses.
\begin{definition}[Support for \inqb]\label{def:support inqb}~
\begin{itemize}
\item $M,s\models p\iff \forall w\in s: V(w,p)=1$
\item $M,s\models\bot\iff s=\emptyset$
\item $M,s\models\phi\land\psi\iff M,s\models\phi$ and $M,s\models\psi$
\item $M,s\models\phi\lori\psi\iff M,s\models\phi$ or $M,s\models\psi$
\item $M,s\models\phi\to\psi\iff\forall t\subseteq s: M,t\models\phi$ implies $M,t\models\psi$
\end{itemize}
\end{definition}
\noindent
The resulting support relation has the following two key features. 
\begin{itemize}
\item Persistence: if $M,s\models\phi$ and $t\subseteq s$, $M,t\models\phi$;
\item Empty State Property: $M,\emptyset\models\phi$ for any formula $\phi$.
\end{itemize}
The notions of logical entailment and equivalence are defined in terms of support in the natural way. More precisely, if $\Phi\subseteq\L$ and $\psi,\chi\in\L$, we let:
\begin{itemize}
\item $\Phi\models\psi\iff\text{for all models $M$ and states $s$: $M,s\models\phi$ for all $\phi\in\Phi$ implies $M,s\models\psi$}$
\item $\psi\equiv\chi\iff \text{for all models $M$ and states $s$: $(M,s\models\psi\iff M,s\models\chi)$}$
\end{itemize}
Although the primitive semantic notion in inquisitive logic is support at an information state, a notion of truth at a world is recovered by defining truth at $w$ as support relative to the corresponding singleton state, $\{w\}$.
\begin{definition}[Truth at a world] $M,w\models\phi\overset{\text{def}}{\iff}M,\{w\}\models\phi$
\end{definition}
\noindent
It is easy to check that, with this definition, all standard connectives obey their usual truth tables; for instance $M,w\models\neg\phi\iff M,w\not\models\phi$, and similarly for $\land,\lor,\to,\leftrightarrow$. As a consequence, all standard propositional formulas receive the usual truth conditions. Moreover, an inquisitive disjunction $\phi\lori\psi$ has the same truth conditions as the declarative disjunction $\phi\lor\psi$, though as we will see, the two usually differ in their support conditions. 

\begin{definition}[Truth-set] The truth-set of a formula $\phi$ in a model $M$, denoted $|\phi|_M$, is the set of worlds in $M$ where $\phi$ is true:
$$|\phi|_M=\{w\in W\mid M,w\models\phi\}$$
\end{definition}

\noindent
For some but not all formulas, support conditions are fully determined by truth conditions in a specific way: support at a state boils down to truth at each world in the state. In this case, we say that the formula is \emph{truth-conditional}.\footnote{In the tradition of dependence logic, truth-conditional formulas are called \emph{flat}.}

\begin{definition}[Truth-conditional formulas] A formula $\phi\in\L$ is \emph{truth-conditional} if for every model $M$ and information state $s$ we have: $M,s\models\phi\iff\forall w\in s: M,w\models\phi$.
\end{definition}

\noindent In inquisitive logic, truth-conditionality is regarded as the key property distinguishing statements (which are truth-conditional) from questions (which are not). 

It will be useful to remark explicitly that truth-conditionality implies (and in fact, given persistence, is equivalent to) a union-closure property.

\begin{observation}[Truth-conditionality implies union-closure]\label{obs:union-closure} Suppose $\phi\in\L$ is truth-conditional. For any model $M$ and any family $S\subseteq\wp(W)$ of information states we have:
$$M,s\models\phi\text{ for all $s\in S$ }\;\Longrightarrow\; M,\bigcup S\models\phi$$
\end{observation}

\noindent 
Standard (i.e., $\lori$-free) propositional formulas in \inqb\ are taken to regiment statements, and indeed, one may verify that they are truth-conditional. Thus, for these formulas the results of our support semantics are fully inter-derivable with those of standard truth-conditional semantics. By contrast, formulas containing $\lori$ are typically not truth-conditional. For an illustration, take the formula $?p$ (shorthand for $p\lori\neg p$). Although this formula is true at every world in every model, its support conditions are non-trivial:
\begin{eqnarray*}
M,s\models{?p}&\iff & M,s\models p\text{ or }M,s\models\neg p\\
&\iff & (\forall w\in s: V(w,p)=1)\text{ or }(\forall w\in s: V(w,p)=0)\\
&\iff & \forall w,w'\in s: V(w,p)=V(w',p)
\end{eqnarray*}
In words, $?p$ is supported by a state $s$ if the truth value of $p$ is \emph{settled} in $s$, i.e., if it is the same at each world $w\in s$.

\subsection{Inquisitive modal logics}

An inquisitive modal logic $\inqml_{\MM}$ extends \inqb\ with a repertoire \MM\ of modal operators $\M$, each with a specified arity $n$. Thus, the language $\L_\MM$ of $\inqml_\MM$ is given by the following definition, where $p\in\P$ and $\M\in\MM$ is an $n$-ary modality:
$$\phi\;:=\;p\mid\bot\mid(\phi\land\phi)\mid(\phi\to\phi)\mid(\phi\lori\phi)\mid \M(\underbrace{\phi,\dots,\phi}_{\text{$n$ times}})$$
In principle, one may consider all sorts of modal operators. However, the focus of inquisitive modal logic is on \emph{declarative} modalities, i.e., modalities which combine with one or more arguments to yield a statement---i.e., a formula which is truth-conditional. 
\begin{definition}[Declarative operators] A $n$-ary operator $\M$ is said to be \emph{declarative} if the formula $\M(\phi_1,\dots,\phi_n)$ is truth-conditional for any $\phi_1,\dots,\phi_n\in\L_\MM$.
\end{definition}

\noindent
Note that in order to define the semantics of a declarative modality $\M$, it suffices to give the truth conditions for formulas $\M(\phi_1,\dots,\phi_n)$; the corresponding support conditions will then simply say that the formula is supported at $s$ if it is true at each world $w\in s$. 

From now on, when we speak of an inquisitive modal logic $\inqml_\MM$ we will assume that $\MM$ is a set of declarative modalities. 
Such a logic is guaranteed to preserve the key properties of support, namely, Persistence and the Empty State Property. 
Moreover, in such a logic we can identify an important syntactic fragment, consisting of formulas that we call \emph{declaratives}.
\begin{definition}[Declaratives]
The set $\L_\MM^!$ of declaratives of an inquisitive modal logic $\inqml_\MM$ is given by the following BNF definition, where $p\in\P$, $\M\in\MM$, and $\phi\in\L_\MM$ stands for an arbitrary formula:
$$\aa\;:=\;p\mid\bot\mid(\aa\land\aa)\mid(\aa\to\aa)\mid\M(\phi,\dots,\phi)$$
\end{definition}
\noindent
More simply put, a formula $\aa$ is a declarative if all the occurrences of $\lori$ in $\aa$, if any, are within the scope of a modality. From now on, we are going to use $\aa,\bb,\gg$ as meta-variables ranging over declaratives, in contrast to $\phi,\psi,\chi$ which range over arbitrary formulas.

The crucial property of the declarative fragment is that, up to equivalence, declaratives represent all and only the statements in the language.

\begin{proposition}\label{prop:declarative variant} Every $\aa\in\L_\MM^!$ is truth-conditional. Moreover, every truth-conditional formula $\phi\in\L_\MM$ is equivalent to some declarative $\phi^!\in\L_\MM^!$.
\end{proposition}

\begin{proof} The first claim is proved by a straightforward induction on the structure of declaratives, using the assumption that modal formulas are always truth-conditional. For the second claim, we associate to each $\phi\in\L_\To$ a declarative $\phi^!$ as follows:
\begin{itemize}
\item $\phi^!=\phi$ if $\phi$ is an atom, $\bot$, or a modal formula $\M(\psi_1,\dots,\psi_n)$
\item $(\phi\land\psi)^!=\phi^!\land\psi^!$
\item $(\phi\lori\psi)^!=\phi^!\lor\psi^!$
\item $(\phi\to\psi)^!=\phi^!\to\psi^!$
\end{itemize}
A straightforward induction reveals that $\phi$ and $\phi^!$ have the same truth conditions. Now if $\phi$ is truth-conditional, then since $\phi^!$ is truth-conditional as well (since it is a declarative), $\phi$ and $\phi^!$ must also have the same support conditions, so we have $\phi\equiv\phi^!$.
\end{proof}

\noindent
Importantly, with each formula $\phi$ of an inquisitive modal logic $\inqml_{\MM}$ we can associate a finite set $\R(\phi)$ of declaratives, called its \emph{resolutions}, in the following way.

\begin{definition}[Resolutions]~ 
\begin{itemize}
\item $\R(\phi)=\{\phi\}$ if $\phi$ is an atom, $\bot$, or a modal formula $\M(\psi_1,\dots,\psi_n)$
\item $\R(\phi\land\psi)=\{\aa\land\bb\mid \aa\in\R(\phi),\bb\in\R(\psi)\}$
\item $\R(\phi\lori\psi)=\R(\phi)\cup\R(\psi)$
\item $\R(\phi\to\psi)=\{\bigwedge_{\aa\in\R(\phi)}(\aa\to f(\aa))\mid f:\R(\phi)\to\R(\psi)\}$
\end{itemize}
\end{definition}

\noindent
The crucial property of resolutions is given by the following normal form result, whose proof is an obvious adaptation of the one for the \inqb\ \citep[cf.\ Prop.\ 3.6.8 in][]{Ciardelli:23book}. 

\begin{proposition}[Normal form]\label{prop:normal form} For every $\phi\in\L_\MM$, if $\R(\phi)=\{\aa_1,\dots,\aa_n\}$ then
$$\phi\equiv\aa_1\lori\dots\lori\aa_n.$$
\end{proposition}
\noindent
Thus, a formula of an inquisitive modal logic (based on declarative modalities) is always equivalent to an inquisitive disjunction of declaratives. Note however that, whereas the resolutions of a propositional formula are always $\lori$-free, the resolutions of a modal formula may contain occurrences of $\lori$, though only within the scope of a modal operator.

\section{Inquisitive neighborhood logic}
\label{sec:inqnl}

Having reviewed some general background on inquisitive modals logics, we may now turn our attention to the specific inquisitive modal logic \inqnl\ that is the topic of this paper.

\paragraph{Syntax.} \inqnl\ is an inquisitive modal logic equipped with a binary declarative modality $\To$, read `yields'. Explicitly, the language $\L_\To$ is thus given by the following definition:
$$\phi\;:=\;p\mid\bot\mid(\phi\land\phi)\mid(\phi\to\phi)\mid(\phi\lori\phi)\mid (\phi\To\phi)$$
Two unary modalities $\ibox$ (`window') and $\idia$ (`kite') are defined in terms of $\To$ as follows:
$$\ibox\phi\;:=\;(\top\To\phi)\qquad\qquad{\idia\phi}\;:=\;\neg (\phi\To\bot)$$
We will also consider some fragments of $\inqnl$, namely, $\inqml_\ibox$, $\inqml_\idia$, and $\inqml_{\ibox\idia}$, based on the following languages:

\begin{center}
\begin{tabular}{l c c c l}
$\L_\ibox$ & $\quad$ & $\phi$ & $:=$ & $p\mid \bot\mid(\phi\land\phi)\mid(\phi\to\phi)\mid(\phi\lori\phi)\mid \ibox\phi$\\[.2cm]
$\L_\idia$ & $\quad$ & $\phi$ & $:=$ & $p\mid \bot\mid(\phi\land\phi)\mid(\phi\to\phi)\mid(\phi\lori\phi)\mid \idia\phi$\\[.2cm]
$\L_{\ibox\idia}$ & $\quad$ & $\phi$ & $:=$ & $p\mid \bot\mid(\phi\land\phi)\mid(\phi\to\phi)\mid(\phi\lori\phi)\mid \ibox\phi\mid \idia\phi$
\end{tabular}
\end{center}

\noindent
The corresponding declarative fragments $\L_\To^!, \L_\ibox^!$, \dots, are defined as in the previous section by restricting to formulas where $\lori$ can occur only within the scope of a modality. 

Each formula $\phi$ of our language is assigned a modal depth $\md(\phi)\in\mathbb{N}$ that tracks the maximum number of nestings of modalities, in the standard way: 
\begin{itemize}
\item $\md(p)=0$ if $p$ is an atom or $\bot$
\item $\md(\phi\circ\psi)=\max(\md(\phi),\md(\psi))$ if $\circ\in\{\land,\lori,\to\}$
\item $\md(\phi\To\psi)=\max(\md(\phi),\md(\psi))+1$
\end{itemize}
The set of formulas of modal depth up to $n$ is denoted by $\L_\To^n$.

\paragraph{Models.}
\inqnl\ is interpreted over neighborhood models. For simplicity, throughout most of the paper, we restrict attention to the case in which neighborhoods are nonempty. This is a modest restriction, as empty neighborhoods do not play a role in any of the standard interpretations of neighborhood models. We will, nevertheless, see in Section \ref{sec:empty} that the approach extends easily to the case in which neighborhoods may be empty. 

\begin{definition}[Inhabited neighborhood models]
An \emph{inhabited neighborhood model}, or \emph{in-model} for short, is a triple $M=\<W,\Sigma,V\>$ where:
\begin{itemize}
\item $W$, the \emph{universe}, is a nonempty set whose elements are called \emph{worlds};
\item $\Sigma:W\to\wp\wp_0(W)$, the \emph{neighborhood map}, assigns to each world $w$ a set $\Sigma(w)$ of neighborhoods, where each neighborhood is a nonempty set of worlds;\footnote{Given a set $X$, we denote by $\wp(X)$ the powerset of $X$, i.e., the set of all subsets of $X$, and we denote by $\wp_0(X)=\wp(X)-\{\emptyset\}$ the set of all nonempty subsets of $X$.}
\item $V:W\times\P\to\{0,1\}$, the \emph{valuation function}, assigns to each world $w$ and atomic sentence $p$ a truth value $V(w,p)$.
\end{itemize}
If $M=\<W,\Sigma,V\>$ is an in-model, the pair $\<W,\Sigma\>$ is called the \emph{frame} underlying $M$.

A \emph{world-pointed model} is a pair $\<M,w\>$ where $M$ is an in-model and $w$ a world in~$M$; similarly, a \emph{state-pointed model} is a pair $\<M,s\>$ where $M$ is an in-model and $s$ is an information state in $M$. 
\end{definition}
\noindent 
Note that an in-model $M=\<W,\Sigma,V\>$ can always be associated with an underlying Kripke model $M_K=\<W,R_{\Sigma},V\>$ obtained by letting the set of successors of $w$ be given by the union of the neighborhoods of $w$: $R_\Sigma[w]=\bigcup\Sigma(w)$.

\paragraph{Semantics.}
Since \inqnl\ is an inquisitive logic, its semantics is given by a recursive definition of support relative to an information state $s\subseteq W$. To the support clauses in Definition \ref{def:support inqb}, we must add a support clause for modal formulas $\phi\To\psi$. Since $\To$ is a declarative modality, such formulas are required to be truth-conditional, and therefore, it suffices to specify their truth conditions relative to a world. These are given by the following clause:
$$M,w\models(\phi\To\psi)\iff \forall s\in\Sigma(w): M,s\models\phi\text{ implies }M,s\models\psi$$
In words, $\phi\To\psi$ 
is true at a world $w$ if every neighborhood of $w$ that supports $\phi$ also supports $\psi$. Spelling out the definitions, one may verify that the operators $\ibox$ and $\idia$ are, respectively, a universal and an existential quantifier over neighborhoods:\footnote{If the empty set were allowed as a neighborhood, the truth conditions for $\idia\phi$ would have to be formulated more carefully as: $M,w\models\idia\phi\iff\text{there is a \emph{nonempty} }s\in\Sigma(w), s\neq\emptyset: M,s\models\phi$. We will come back to this point in Section \ref{sec:empty}.}
\begin{itemize}
\item $M,w\models\ibox\phi\iff\forall s\in\Sigma(w): M,s\models\phi$
\item $M,w\models\idia\phi\iff\exists s\in\Sigma(w): M,s\models\phi$
\end{itemize}
In the rest of this section, we familiarize ourselves with the workings of this logic by considering a range of examples and special cases. We start by examining the fragments $\inqml_\idia$ and $\inqml_\ibox$ containing only one of the modalities $\idia$ and $\ibox$.

\paragraph{$\inqml_\idia$.} Considering first a formula $\idia\aa$ where $\aa$ is a declarative. Using the fact that $\aa$ is truth-conditional, we have:
\begin{eqnarray*}
M,w\models\idia\aa&\iff & \exists s\in\Sigma(w): M,s\models\aa\\
&\iff & \exists s\in\Sigma(w)\text{ s.t.\ }\forall w\in s: M,w\models\aa
\end{eqnarray*}
In words, $\idia\aa$ is true at $w$ in case $\aa$ is true throughout some neighborhood of $w$. This recovers the $\exists\forall$-interpretation of the modal formula $\Box\aa$ discussed in the introduction. Thus, simply by replacing the modality $\kbox$ by $\idia$, we may identify the standard modal language $\L_\Box$ equipped with the $\exists\forall$-semantics with the $\lori$-free fragment of $\L_{\idia}$. 

Conversely, we can show that any declarative in  $\inqml_{\idia}$ is equivalent to a formula of $\L_\kbox$ under the $\exists\forall$-semantics.  
To see why, start with the following proposition. 

\begin{proposition}\label{prop:idia dist} For any $\phi,\psi\in\L_\To$, $\idia(\phi\lori\psi)\equiv\idia\phi\lor\idia\psi$.
\end{proposition} 

\begin{proof}
Since the relevant formulas are truth-conditional, it suffices to check that they have the same truth conditions. For any model $M$ and world $w$, we have:
\begin{eqnarray*}
M,w\models\idia(\phi\lori\psi)&\iff& \exists s\in\Sigma(w): M,s\models\phi\lori\psi\\
&\iff& \exists s\in\Sigma(w): M,s\models\phi\text{ or }M,s\models\psi\\
&\iff& (\exists s\in\Sigma(w): M,s\models\phi)\text{ or }(\exists s\in\Sigma(w):  M,s\models\psi)\\
&\iff& M,w\models\idia\phi\text{ or }M,w\models\idia\psi\\
&\iff& M,w\models\idia\phi\lor\idia\psi
\end{eqnarray*}
\end{proof}

\noindent
From the above equivalence and the normal form result in Proposition \ref{prop:normal form} we then get the following corollary, which ensures that inquisitive disjunctions occurring under $\idia$ can always be paraphrased away. 

\begin{corollary}\label{cor:idia elimination} Every declarative in $\inqml_{\idia}$ is equivalent to a $\lori$-free one.
\end{corollary}

\begin{proof} Let $\phi$ be a formula and let $\R(\phi)=\{\aa_1,\dots,\aa_n\}$. Combining propositions \ref{prop:normal form} and \ref{prop:idia dist} we obtain:
\begin{eqnarray*}
\idia\phi&\equiv&\idia(\aa_1\lori\dots\lori\aa_n)\\
&\equiv&\idia\aa_1\lor\dots\lor\idia\aa_n
\end{eqnarray*}
Now take an arbitrary declarative $\aa\in\L_{\idia}^!$. Applying the above equivalence recursively, from the inner-most occurrence of $\idia$ in $\aa$ to the outer-most one, one can ultimately remove every occurrence of $\lori$ in $\aa$. 
\end{proof}

\noindent
Thus, up to logical equivalence, the declarative fragment of $\inqml_\idia$ can be identified with the standard modal language $\L_\Box$ equipped with the $\exists\forall$-semantics.

\paragraph{$\inqml_\ibox$.} Next, consider the semantics of a formula $\ibox\aa$, where $\aa$ is a declarative. Using the truth-conditionality of $\aa$, we have:
\begin{eqnarray*}
M,w\models\ibox\aa&\iff&\forall s\in \Sigma(w): M,s\models\aa\\
&\iff&\forall s\in\Sigma(w)\forall v\in s: M,v\models\aa\\
&\iff&\forall v\in \bigcup\Sigma(w): M,v\models\aa\\
&\iff&\forall v\in R_\Sigma[w]: M,v\models\aa
\end{eqnarray*}
In words, $\ibox\aa$ is true at $w$ if $\aa$ is true at every successor of $w$ in the underlying Kripke model. In other words, when applied to declaratives, $\ibox$ behaves like a standard $\Box$ modality interpreted over the underlying Kripke model $M_K$. 

Things become more interesting when $\ibox$ is applied to a question. For an illustration, consider the formula $\ibox{?p}$. Recalling that $?p$ is supported by a state $s$ just in case $p$ has the same truth value in all worlds in $s$, we have:
\begin{eqnarray*}
M,w\models\ibox{?p} &\iff&\forall s\in \Sigma(w): M,s\models{?p}\\
&\iff&\forall s\in\Sigma(w):\text{$p$ has the same truth value in all $w\in s$}
\end{eqnarray*}
In words, $\ibox?p$ says that the truth value of $p$ is constant in each neighborhood, i.e., each neighborhood ``settles'' the truth value of $p$. Thus, for instance, in Fig.~\ref{fig:1}, the formula $\ibox?p$ is true at $w_1$, but false at $w_2$ and $w_3$. 

This example already suffices to show that, unlike in the case of $\idia$, modal claims $\ibox\phi$ involving questions do not in general reduce to Boolean combinations of modal claims $\ibox\aa$ involving declaratives. For as we have seen, the latter claims only depend for their truth at $w$ on the set of successors $R_\Sigma[w]$. Since the three worlds  in Fig.~\ref{fig:1} have the same set of successors (the union of the neighborhoods is the same in each case), these worlds will satisfy the same formulas of the form $\ibox\aa$ where $\aa$ is a declarative. Yet, as we saw, $w_1$ and $w_2$ are distinguished by the formula $\ibox?p$. Therefore, this formula cannot be equivalent to any Boolean combination of formulas $\ibox\aa$ where $\aa$ is declarative. 

\paragraph{General case.} 
Consider a modal formula $(\aa\To\bb)$ where $\aa$ and $\bb$ are declaratives. This formula expresses a kind of global consequence from the perspective of the world $w$: for every neighborhood $s\in\Sigma(w)$, if $\aa$ is true everywhere in $s$, $\bb$ is true everywhere in $s$.
$$M,w\models (\aa\To\bb)\iff \forall s\in\Sigma(w): s\subseteq|\aa|_M\text{ implies }s\subseteq|\bb|_M$$
More generally, if $\aa,\bb_1,\dots,\bb_n$ are declaratives, $(\aa\To \bb_1\lori\dots\lori\bb_n)$ expresses a kind of multiple-conclusion global consequence: for every neighborhood $s$, if $\aa$ is true everywhere in $s$, some $\bb_i$ is true everywhere in $s$.
$$M,w\models(\aa\To\bb_1\lori\dots\lori\bb_n)
\;\iff\; \forall s\in \Sigma(w): s\subseteq|\aa|_M\text{ implies }(s\subseteq|\bb_i|_M\text{ for some }i)$$
An example of this kind of formula is $p\To {?q}$ (short for $p\To q\lori\neg q$), which allows us to distinguish the worlds $w_2$ and $w_3$ in Fig.~\ref{fig:1}: in $w_2$, any neighborhood where $p$ is everywhere true is a singleton, and must thus settle the truth value of $q$; thus, the formula is true in $w_2$. However, it is not true in $w_3$, since $p$ is everywhere true in the neighborhood $\{w_{pq},w_{p\overline q}\}$, which does not settle the truth value of $q$. 

It is also interesting to consider the semantics of the formula $\neg(\aa\To(\neg\bb_1\lori\dots\lori\neg\bb_n))$, where again $\aa,\bb_1,\dots,\bb_n$ are declaratives. We have:
\begin{eqnarray*}
M,w\models\neg(\aa\To\Lori_{i\le n}\neg\bb_i)&\iff &\exists s\in\Sigma(w): M,s\models\aa\text{ and }M,s\not\models\neg\bb_i\text{ for each }i\\
&\iff &\exists s\in\Sigma(w): s\subseteq|\aa|_M\text{ and }s\not\subseteq|\neg\bb_i|_M\text{ for each }i\\
&\iff &\exists s\in\Sigma(w): s\subseteq|\aa|_M\text{ and }s\cap|\bb_i|_M\neq\emptyset\text{ for each }i
\end{eqnarray*}
Thus, the given formula expresses the existence of a neighborhood $s$ such that $\aa$ is true everywhere in $s$ and for each $i\le n$, $\bb_i$ is true somewhere in $s$. This is exactly the meaning of a modal formula of the form $\Box(\beta_1,\dots,\beta_n;\aa)$ in the instantial neighborhood logic of \cite{Benthem:17}; we will come back to the connection in Section 
\ref{sec:inl}.

In contrast to inquisitive disjunctions on the right of $\To$, which play an essential role in our logic (cf.\ the proof of Corollary \ref{cor:expressivity gap} in the next section), inquisitive disjunctions on the left of $\To$ can always be paraphrased away, as the following propositions shows.

\begin{proposition} For any $\phi,\psi,\chi\in\L_\To$ we have $(\phi\lori\psi\To\chi)\equiv(\phi\To\chi)\land(\psi\To\chi)$.
\end{proposition}

\begin{proof} Immediate by spelling out the truth conditions of the relevant formulas.
\end{proof}

\noindent
Using this equivalence in combination with the normal form result in Proposition \ref{prop:normal form}, we obtain the following corollary (cf.\ the proof of Corollary \ref{cor:idia elimination}).

\begin{corollary} Any formula of \inqnl\ is equivalent to one where $\lori$ does not occur on the left of $\To$.
\end{corollary}

\noindent
An interesting feature of the logic $\inqnl$ is the interplay of the two conditionals $\To$ and $\to$, that can be used to specify global and local restrictions respectively: while $\To$ allows us to restrict the set of neighborhoods under consideration, $\to$ allows us to restrict the worlds in each individual neighborhood. For an example, consider the formula $(p\To (q\to{?r}))$. This expresses the following fact: if we restrict to neighborhoods where $p$ is everywhere true, and then we restrict each of these neighborhoods to the $q$-worlds, all the resulting states  settle whether~$r$:
$$M,w\models (p\To (q\to{?r}))\iff \forall s\in\Sigma(w): M,s\models p\text{ implies }M,s\cap|q|_M\models{?r}$$
Similarly, the formula $({?p}\To ({?q}\to{?r}))$ expresses that if we restrict to neighborhoods that settle whether $p$, and then look at the parts of such neighborhoods where the truth value of $q$ is settled, each of these parts settles whether~$r$:
\begin{eqnarray*}
M,w\models ({?p}\To ({?q}\to{?r}))&\iff &\forall s\in\Sigma(w): M,s\models {?p} \text{ implies }\\
&&\phantom{\forall s\in\Sigma(w):\;}M, s\cap|q|_M\models{?r}\,\text{ and }\,M,s\cap|\neg q|_M\models{?r}
\end{eqnarray*}
Finally, note that in \inqnl\ we may define a Lewisian counterfactual operator. Indeed, if $\aa$ and $\bb$ are declaratives, we may define the counterfactual $\aa\,\counterfactual\,\bb$ as follows:
$$\aa\,\counterfactual\,\bb\;:=\;\ibox\neg\aa\lor\neg((\aa\to\bb)\To\neg\aa)$$
It is an easy exercise to verify that this formula receives precisely the truth conditions that Lewis assigns to a counterfactual, namely: 
\begin{eqnarray*}
M,w\models \aa\,\counterfactual\,\bb &\iff & \text{$\forall s\in\Sigma(w): \; s\cap|\aa|_M=\emptyset$, or}\\
& & \text{$\exists s\in\Sigma(w): (s\cap|\aa|_M\neq\emptyset$ and $s\cap|\aa|_M\subseteq|\bb|_M$})
\end{eqnarray*}

\section{Monotonicity properties of fragments} 
\label{sec:fragments}

A couple of questions may arise naturally at this point. Is it really necessary to base our language on the binary modality $\To$? Could we, instead, take the unary modalities $\ibox$ and $\idia$ as primitives? Also, is it crucial that the modality $\To$ is added to an underlying inquisitive logic? Could we, instead, just add $\To$ to standard propositional logic? 

In this section, we will answer these questions by showing that various syntactic fragments of the language $\L_\To$ are characterized by salient monotonicity properties with respect to the set of neighborhoods---properties which are not shared by the full language. The relevant fragments are the languages $\L_\ibox,\L_\idia,\L_{\ibox\idia}$ defined in the previous section, as well as the fragment $\L_{\To!}$ where only declaratives are allowed to occur on the right of~$\To$ (i.e., $\L_{\To!}$ includes only modal formulas of the form $(\phi\To\aa)$ where $\aa$ is a declarative).

To define the relevant monotonicity properties, we first introduce some closure operations on the set of neighborhoods. 

\begin{definition}[Monotonic closures]~\\
Given an in-model $M=\<W,\Sigma,V\>$, we define its:
\begin{itemize}
\item upward closure $M^\uparrow=\<W,\Sigma^\uparrow,V\>$,\\ where
$\Sigma^\uparrow(w)=\{t\in\wp_0(W)\mid t\supseteq s\text{ for some }s\in\Sigma(w)\}$;
\item downward closure $M^\downarrow=\<W,\Sigma^\downarrow,V\>$,\\ where 
$\Sigma^\downarrow(w)=\{t\in\wp_0(W)\mid t\subseteq s\text{ for some }s\in\Sigma(w)\}$;
\item convex closure $M^\updownarrow=\<W,\Sigma^\updownarrow,V\>$,\\ where 
$\Sigma^\updownarrow(w)=\{t\in\wp_0(W)\mid s\subseteq t\subseteq s'\text{ for some }s,s'\in\Sigma(w)\}$;
\item union closure $M^\cup=\<W,\Sigma^\cup,V\>$,\\ where
$\Sigma^\cup(w)=\{t\in\wp_0(W)\mid t=\bigcup S\text{ for some }S\subseteq\Sigma(w)\}$.

\end{itemize}
We say that a model $M$ is \emph{upward-monotone} if $M=M^\uparrow$, \emph{downward-monotone} if $M=M^\downarrow$, \emph{convex} if $M=M^\updownarrow$, and \emph{union-closed} if $M=M^\cup$.
\end{definition}

\begin{definition}[Monotonicity properties] We say that a formula $\phi\in\L_{\To}$ is:
\begin{itemize}
\item upward-invariant if for any $M,s$ we have $M,s\models\phi\iff M^\uparrow,s\models\phi$
\item downward-invariant if for any $M,s$ we have $M,s\models\phi\iff M^\downarrow,s\models\phi$
\item convex-invariant if for any $M,s$ we have $M,s\models\phi\iff M^\updownarrow,s\models\phi$
\item union-invariant if for any $M,s$ we have $M,s\models\phi\iff M^\cup,s\models\phi$
\end{itemize}
\end{definition}

\noindent
Note that upward- or downward-invariance imply convex-invariance. For instance, suppose that $\phi$ is upward-invariant; since $M^\uparrow=(M^\updownarrow)^\uparrow$, for any $M,s$ we have:
$$M,s\models\phi\iff M^{\uparrow},s\models\phi\iff (M^{\updownarrow})^{\uparrow},s\models\phi\iff M^{\updownarrow},s\models\phi$$
Similarly, upward-invariance implies union-invariance, since $M^\uparrow=(M^\cup)^\uparrow$. 

We are now ready to state the monotonicity properties that characterize the fragments of \inqnl\ we considered.

\begin{proposition}[Monotonicity properties of fragments]\label{prop:monotonicity}~
\begin{itemize}
\item All formulas in $\L_{\idia}$ are upward-invariant (thus, also convex- and union-invariant).
\item All formulas in $\L_\ibox$ are downward-invariant (thus, also convex-invariant).
\item All formulas in $\L_{\ibox\idia}$ are convex-invariant.
\item All formulas in $\L_{\To!}$ are union-invariant.
\end{itemize}
\end{proposition}

\begin{proof} 
We prove only the first and the last claim; the proofs of the second and the third claim are analogous to the one given for the first. 

To establish the first claim, we proceed by induction on $\phi\in\L_{\idia}$. The only non-trivial case is the inductive step for $\phi=\idia\psi$. Since $\idia\psi$ is truth-conditional, it suffices to show that for all worlds $w$ in $M$ we have $M,w\models\idia\psi\iff M^\uparrow,w\models\idia\psi$.

The $\Rightarrow$ direction of the claim is immediate from the induction hypothesis and the fact that $\Sigma(w)\subseteq\Sigma^\uparrow(w)$. For the $\Leftarrow$ direction, suppose
 $M^{\uparrow},w\models\idia\psi$. Then for some $s\in\Sigma^{\uparrow}(w)$ we have $M^{\uparrow},s\models\psi$. By induction hypothesis, $M,s\models\psi$. By definition of $\Sigma^\uparrow(w)$, we have $s'\subseteq s$ for some $s'\in\Sigma(w)$, and by the persistence of support, $M,s'\models\psi$. Therefore, $M,w\models\idia\psi$.

For the last claim, again we proceed by induction on $\phi\in\L_{\To!}$. The only non-trivial case is the inductive step for a modal formula $\phi=(\psi\To\aa)$, where by definition of $\L_{\To!}$, $\aa$ must be a declarative and thus truth-conditional by Proposition \ref{prop:declarative variant}. 
Since the formula $(\psi\To\aa)$ is truth-conditional, it suffices to show that for all worlds $w$ in $M$ we have $M,w\models(\psi\To\aa)\iff M^\cup,w\models(\psi\To\aa)$.

The $\Leftarrow$ direction of the claim is immediate from the induction hypothesis and the fact that $\Sigma(w)\subseteq\Sigma^\cup(w)$. For the $\Rightarrow$ direction, suppose $M,w\models(\psi\To\aa)$. This means that for all $s\in\Sigma(w)$, $M,s\models\psi$ implies $M,s\models\aa$. By the induction hypothesis, this means that for all $s\in\Sigma(w)$, $M^\cup,s\models\psi$ implies $M^\cup,s\models\aa$ (call this fact $(*)$). Now take an arbitrary $t\in\Sigma^\cup(w)$: by definition of $\Sigma^\cup(w)$, $t=\bigcup S$ for some $S\subseteq \Sigma(w)$. Now we have:
$$
\begin{array}{rcll}
M^\cup,t\models\psi &\Longrightarrow& \forall s\in S: M^\cup,s\models\psi & \text{by persistence, since $t=\bigcup S$}\\ 
&\Longrightarrow& \forall s\in S: M^\cup,s\models\aa & \text{from  $(*)$ and the fact that $S\subseteq\Sigma(w)$}\\ 
&\Longrightarrow& M^\cup,t\models\aa & \text{by Observation \ref{obs:union-closure} since $t=\bigcup S$} 
\end{array}
$$
This shows that $M^\cup,w\models(\psi\To\aa)$, as desired.

\end{proof}

\noindent
In words, the proposition says that $\idia$-formulas are insensitive to the addition of supersets of neighborhoods, $\ibox$-formulas are insensitive to the addition of subsets of neighborhoods, formulas containing both $\ibox$ and $\idia$ are insensitive to the addition of states which are included in between two neighborhoods, and formulas where questions do not appear on the right of $\To$ are insensitive to the addition of states which are unions of neighborhoods.

As a corollary, we also get that the full language $\L_\To$ is strictly more expressive than any of the fragments we considered.

\begin{corollary}\label{cor:expressivity gap} $\L_\To$ is strictly more expressive than any of $\L_{\idia}$, $\L_{\ibox}$, $\L_{\ibox\idia}$, $\L_{\To!}$.
\end{corollary}

\smallskip\noindent\emph{Proof.}
It suffices to note that $\L_\To$ as a whole does not have any of the monotonicity properties that characterize the sub-fragments. This is illustrated by the following examples.
\begin{itemize}
\item The formula $(p\To q)$ is union-invariant by Proposition \ref{prop:monotonicity}, but not it is easy to see that it is not convex-invariant (and, thus, also not upward- or downward-invariant). 

\item The formula $\ibox{?p}$ is downward-invariant by Proposition \ref{prop:monotonicity}, but it is easy to see that it is not union-invariant (and, thus, not upward-invariant). 
Note that this means that any formula equivalent to $\ibox{?p}$ must contain an occurrence of a question on the right of $\To$. 
\hfill$\Box$
\end{itemize}

\noindent This result shows that, to get the full expressive power of our logic, it is crucial to use the binary modality $\To$ in combination with questions. Having only the unary modalities $\ibox$ and $\idia$ is not enough, nor is having the binary modality $\To$ in the absence of questions. 

As mentioned in the introduction, most previous work on inquisitive modal logic has focused on a language based on the modality $\ibox$. However, the semantics was restricted to \emph{downward-monotone} neighborhood models. In that context, $\To$ becomes definable in terms of $\ibox$, as the following proposition shows, and therefore, the logics $\inqml_\ibox$ and $\inqnl$ are equi-expressive.

\begin{proposition}[$\To$ is definable from $\ibox$ over downward-monotone models] Let $M$ be a downward-monotone model and $w$ a world in $M$. For any formulas $\phi,\psi\in\L_\To$, we have:
$$M,w\models(\phi\To\psi)\iff M,w\models\ibox(\phi\to\psi)$$
\end{proposition}

\begin{proof} Spelling out the truth conditions of $\ibox(\phi\to\psi)$, we have:
$$M,w\models\ibox(\phi\to\psi)\iff\forall s\in\Sigma(w)\forall t\subseteq s:M,t\models\phi\text{ implies }M,t\models\psi$$
The second quantification can be restricted to non-empty subsets of $s$, since the empty state trivially supports any formula. However, in the context of a downward-monotone model, every non-empty subset of a neighborhood is itself a neighborhood. Thus, the above condition can be simplified to ``$\forall s\in\Sigma(w):M,s\models\phi\text{ implies }M,s\models\psi$''
which is precisely what is required to have $M,w\models(\phi\To\psi)$
\end{proof}

\noindent
Interestingly, analogous definitions of $\To$ in terms of $\idia$ over upward-monotone models, or in terms of both $\ibox$ and $\idia$ over convex models, do not seem to be available; we leave this as a conjecture to be explored in future work.

\section{Bisimilarity and expressive power}
\label{sec:expressive power}

In this section, we introduce the natural notion of bisimilarity for neighborhood models \citep[cf.\ also][]{Benthem:17,CiardelliOtto:18}, and characterize the expressive power of \inqnl\ in terms of it. 

\paragraph{Bisimilarity.} The notion of bisimilarity captures the idea of \emph{behavioral equivalence}. In the setting of neighborhood models, it is natural to regard two worlds as behaviorally equivalent when they agree on atomic propositions and every neighborhood of the one is behaviorally equivalent to a neighborhood of the other. In turn, it is natural to regard two neighborhoods as behaviorally equivalent when each world in the one is behaviorally equivalent to a world in the other. This leads naturally to a notion of bisimilarity, which may be formalized equivalently in two ways: by defining a notion of \emph{bisimulation relation}, or in terms of a \emph{bisimulation game}. We start with the former approach.

\begin{definition}[Egli-Milner lifting] If $R\subseteq X\times Y$ is a relation, its Egle-Miler lifting is the relation $\overline R=\wp(X)\times\wp(Y)$ defined by letting, for every $A\subseteq X$ and $B\subseteq Y$:
\begin{eqnarray*}
A\overline R B&\iff&\text{for every } a\in A\text{ there is some } b\in B\text{ with }aRb\text{ and }\\
&&\text{for every } b\in B\text{ there is some } a\in A\text{ with }aRb
\end{eqnarray*}
\end{definition}

\begin{definition}[Bisimulation]
Let $M=\<W,\Sigma,V\>$ and $M'=\<W',\Sigma',V'\>$ be two in-models. A bisimulation between $M$ and $M'$ is a non-empty relation $Z\subseteq W\times W'$ such that whenever $wZw'$ holds, the following conditions are satisfied:
\begin{itemize}
\item Atomic equivalence: for each $p\in\P$, $V(w,p)=V'(w',p)$;
\item Forth condition: for every $s\in\Sigma(w)$ there is some $s'\in\Sigma'(w')$ such that $s\overline Zs'$;
\item Back condition: for every $s'\in\Sigma'(w')$ there is some $s\in\Sigma(w)$ such that $s\overline Zs'$.
\end{itemize}
We say that two world-pointed models $\<M,w\>$ and $\<M',w'\>$ are \emph{bisimilar} (notation: $M,w\sim M',w'$) if there is a bisimulation $Z$ between $M$ and $M'$ with $wZw'$. We say that two state-pointed models $\<M,s\>$ and $\<M',s'\>$ are bisimilar ($M,s\sim M',s'$) if there is a bisimulation $Z$ between $M$ and $M'$ such that $s\overline Z s'$. 
\end{definition}

\noindent
From the definition it follows immediately that $M,w\sim M',w'$ implies the following: 
\begin{itemize}
\item Atomic equivalence: for each $p\in\P$, $V(w,p)=V'(w',p)$;
\item Forth condition: for each $s\in\Sigma(w)$ there is some $s'\in\Sigma'(w')$ such that $M,s\sim M',s'$;
\item Back condition: for each $s'\in\Sigma'(w')$ there is some $s\in\Sigma(w)$ such that $M,s\sim M',s'$.
\end{itemize}
Moreover, it is easy to verify that $M,s\sim M',s'$ holds iff for every $w\in s$ there is a $w'\in s'$ with $M,w\sim M',w'$ and vice versa.

It is also important to consider the $n$-step approximations of bisimilarity, which capture the idea of behavioral equivalence up to $n$ modal transitions (a modal transition being the step from a world to one of its neighborhoods).

\begin{definition}[$n$-Bisimulation]
An $n$-bisimulation between $M$ and $M'$ is a family $(Z_i)_{i\le n}$ of relations $Z_i\subseteq W\times W'$ such that whenever $wZ_iw'$ holds, the following conditions are satisfied:
\begin{itemize}
\item Atomic equivalence: for each $p\in\P$, $V(w,p)=V'(w',p)$;
\item Forth condition: if $i>0$, for each $s\in\Sigma(w)$ there is some $s'\in\Sigma'(w')$ s.t.\ $s\overline Z_{i-1}s'$;
\item Back condition: if $i>0$, for each $s'\in\Sigma'(w')$ there is some $s\in\Sigma(w)$ s.t.\ $s\overline Z_{i-1}s'$.
\end{itemize}
We say that two world-pointed models $\<M,w\>$ and $\<M',w'\>$ are \emph{$n$-bisimilar}, denoted $M,w\sim_n M',w'$, if there is an $n$-bisimulation $(Z_i)_{i\le n}$ between $M$ and $M'$ with $wZ_n w'$. We say that two state-pointed models $\<M,s\>$ and $\<M',s'\>$ are $n$-bisimilar ($M,s\sim_n M',s'$) if  
there is an $n$-bisimulation $(Z_i)_{i\le n}$ with $s\overline Z_n s'$. 
\end{definition}
\noindent
From the definition it follows that $M,w\sim_n M',w'$ implies the following conditions: 
\begin{itemize}
\item Atomic equivalence: for each $p\in\P$, $V(w,p)=V'(w',p)$;
\item Forth condition: if $n>0$, for each $s\in\Sigma(w)$ there is  $s'\in\Sigma'(w')$ s.t.\ ${M,s\sim_{n-1}M',s'}$;
\item Back condition: if $n>0$, for each $s'\in\Sigma'(w')$ there is $s\in\Sigma(w)$ s.t.\ $M,s\sim_{n-1}M',s'$.
\end{itemize}
Moreover, it is easy to check that $M,s\sim_n M',s'$ holds iff for every $w\in s$ there is a $w'\in s'$ with $M,w\sim_n M',w'$ and vice versa.

Bisimilarity can also be given an alternative characterization in terms of a bisimulation game with two players: Spoiler and Duplicator. The game alternates between world-positions $\<w,w'\>\in W\times W'$ and state-positions $\<s,s'\>\in\wp(W)\times\wp(W')$. Playing from a world-position, Spoiler picks a neighborhood of either world and Duplicator responds with a neighborhood of the other world, leading to a state position. Playing from a state position, Spoiler picks a world in either state and Duplicator responds with a world in the other state. A \emph{round} of the game is a sequence of four moves leading from a world-position to the next world-position. If Duplicator is unable to make a move, or if a world-position is reached where the worlds disagree on some atomic sentence, the match ends with a win for Spoiler; if such a situation never occurs in the match, Duplicator wins. We consider two versions of the game: in the $n$-round version, the match ends when $n$ full rounds have been played; in the unbounded version, the match is allowed to go on indefinitely. 

Now we can say that two world-pointed models $\<M,w\>$ and $\<M',w'\>$ 
are bisimilar (respectively, $n$-bisimilar) if Duplicator has a winning strategy in the unbounded game (respectively, $n$-round game) starting from the position $\<w,w'\>$. Similarly, two state-pointed models $\<M,s\>$ and $\<M',s'\>$ are bisimilar ($n$-bisimilar) if Duplicator has a winning strategy in the unbounded ($n$-round) game starting from $\<s,s'\>$.
We leave it to the reader to verify that this characterization of ($n$-)bisimilarity coincides with the one given earlier in terms of ($n$-)bisimulation relations.

As in the case of relational semantics, bisimilarity implies $n$-bisimilarity for each $n\in\mathbb{N}$, but not conversely: two worlds can be $n$-bisimilar for each $n\in\mathbb{N}$, and yet not bisimilar. A straightforward adaptation of the standard example from Kripke semantics suffices to make the point.

\begin{example}[Bisimilarity is strictly stronger than $n$-bisimilarity for each $n$] \label{ex:infinity}
Consider the in-model $M=\<W,\Sigma,V\>$ defined as follows:
\begin{itemize}
\item $W=\{w_0,w_1\}\cup\{v_{ij}\mid 0\le i< j\le\omega\}$
\item $\Sigma(w_0)=\{\{v_{ij}\}\mid 0\le i<j<\omega\}$
\item $\Sigma(w_1)=\{\{v_{ij}\}\mid 0\le i<j\le\omega\}$
\item $\Sigma(v_{ij})=\left\{
\begin{array}{ll}
\{v_{(i+1)j}\} & \text{if $i+1<j$}\\
\emptyset & \text{if $i+1=j$}\\
\end{array}
\right.$
\item $V(p,w)=0$ for every $p\in\P$ and $w\in W$
\end{itemize}
We leave it to the reader to verify that $M,w_0\sim_n M,w_1$ for each $n\in\mathbb{N}$ but $M,w_0\not\sim M,w_1$.
\end{example}

\paragraph{Expressive power of \inqnl.}
We are now going to show that the expressive power of \inqnl\ can be characterized naturally in terms of the notion of bisimilarity.

\begin{definition}[Modal equivalence]
Two world-pointed models $\<M,w\>$ and $\<M',w'\>$ are modally equivalent (notation: $M,w\modeq M',w'$) if they make true the same formulas of \inqnl, and $n$-modally equivalent ($M,w\modeq_n M',w'$) if they make true the same formulas of modal depth up to $n$. Similarly for state-pointed models $\<M,s\>$ and $\<M',s'\>$, where ``make true'' is replaced by ``support''.
\end{definition}

\noindent
The first result we will show is that, provided the set $\P$ of atoms is finite, $n$-bisimilarity coincides precisely with $n$-modal equivalence---that is, the worlds (or states) that can be distinguished by means of formulas of modal depth up to $n$ are precisely those that exhibit a behavioral difference within $n$ steps. 

One direction of the result, which does not require the finiteness of \P, can be proved by a simple induction; we omit the straightforward proof.

\begin{proposition}[$n$-bisimilarity implies $n$-modal equivalence]\label{prop:n-bis to n-modeq}~
\begin{itemize}
\item $M,w\sim_n M',w'$ implies $M,w\modeq_n M',w'$
\item $M,s\sim_n M',s'$ implies $M,s\modeq_n M',s'$.

\end{itemize}
\end{proposition}

\noindent Since full bisimilarity implies $n$-bisimilarity for each $n$, as an immediate corollary we get that \inqnl\ as a whole is bisimulation invariant.

\begin{corollary}[Bisimulation invariance]~
\begin{itemize}
\item $M,w\sim M',w'$ implies $M,w\modeq M',w'$
\item $M,s\sim M',s'$ implies $M,s\modeq M',s'$
\end{itemize}
\end{corollary}

\noindent To show that $n$-modal equivalence implies $n$-bisimilarity if the set \P\ of atoms is finite, we define characteristic formulas $\chi_{M,w}^n$ that capture the equivalence class of the pointed model $\<M,w\>$ modulo $n$-bisimilarity.

\begin{definition}[Characteristic formulas for worlds]~\\
Assume \P\ is finite. We define by simultaneous recursion formulas $\chi_{M,w}^n$ for world-pointed models $\<M,w\>$, and $\pi_{M,s}^n$ for state-pointed models $\<M,s\>$, in the following way:
\begin{itemize}
\item $\chi_{M,w}^0=\bigwedge_{p\in\P}\delta_w(p)$ \hfill where $\delta_w(p)=\left\{
\begin{array}{ll}
p & \text{if $V(w,p)=1$}\\
\neg p & \text{if $V(w,p)=0$}
\end{array}
\right.$
\item $\pi_{M,s}^{n}=\neg(\bigvee\{\chi_{M,w}^n\mid w\in s\}\To \Lori\{\neg\chi_{M,w}^n\mid w\in s\})$
\item $\chi_{M,w}^{n+1}=\chi_{M,w}^0\land \bigwedge\{\pi_{M,s}^n\mid s\in\Sigma(w)\}\land \bigwedge\{\neg\pi_{M',s'}^n\mid\neg\exists s\in\Sigma(w)\text{ with }M,s\sim_n M',s'\}$
\end{itemize}
\end{definition}
\noindent
To verify that the relevant formulas are well-defined, we need to make sure that all conjunctions and disjunctions that appear in the definition are finite.\footnote{Of course, strictly speaking this only defines the relevant formulas up to equivalence. If needed, one may assume that an arbitrary ordering of the formulas in $\L_\To$ is fixed in the background, which determines the order in which the relevant conjuncts or disjuncts are listed.} For this, the crucial observation is that for each number $n$, the set of characteristic formulas $X_n=\{\chi_{M,w}^n\mid\<M,w\>\text{ a world-pointed model}\}$ and $\Pi_n=\{\pi_{M,s}^n\mid\<M,s\>\text{ a state-pointed model}\}$ are finite. For the case $n=0$, the finiteness of $X_0$ is guaranteed by the assumption that \P\ is finite and the fact that inquisitive propositional logic is locally tabular \citep{Ciardelli:23book}. Moreover, the finiteness of a given $X_n$ implies the finiteness of $\Pi_n$, since each formula $\pi\in\Pi_n$ is determined by a set of formulas $Y\subseteq X_n$ via the mapping $Y\mapsto \neg(\bigvee Y\To\Lori \{\neg\chi\mid\chi\in Y\})$, and there are only finitely many distinct subsets of the finite set $X_n$. Finally, for the inductive step, suppose $X_n$ and $\Pi_n$ are finite. Then, the formulas $\chi_{M,w}^{n+1}$ are well-defined, since the relevant conjunctions are over finite sets. Moreover, the set $X_{n+1}$ of such formulas is also finite, since an element of $X_{n+1}$ is a conjunction of one formula $\chi\in X_0$ and a number of distinct formulas in $\Pi_n$ or negations of such formulas; since $X_0$ and $\Pi_n$ are finite, the number of such conjunctions is finite. 

Thus, the formulas $\chi_{M,w}^n$ and $\pi_{M,s}^n$ are indeed well-defined. Note moreover that these formulas are declaratives, as inquisitive disjunction occurs in them only in the consequent of $\To$, and that $\chi_{M,w}^n$ has modal depth $n$. The following lemma states the key property of these formulas.

\begin{lemma}\label{lemma:characteristic formulas} If $\P$ if finite, for every $n\in\mathbb{N}$ we have:
\begin{enumerate}
\item $M,w\models\chi_{M',w'}^n\iff M,w\sim_n M',w'$
\item $M,w\models\pi_{M',s'}^n\iff\text{for some }s\in\Sigma(w): M,s\sim_n M',s'$ 
\end{enumerate}
\end{lemma}

\noindent\emph{Proof.}
We first prove that for a fixed $n$, claim (1) implies claim (2). So, suppose (1) holds for $n$. We prove the two directions of the biconditional in (2).

\begin{itemize}
\item[$\Rightarrow$] Suppose $M,w\models\pi_{M',s'}^n$. By definition of $\pi_{M',s'}^n$, this means that there exists an $s\in\Sigma(w)$ such that $M,s\models\bigvee\{\chi_{M',w'}^n\mid w'\in s'\}$, but $M,s\not\models\Lori\{\neg\chi_{M',w'}^n\mid w'\in s'\}$.

We claim that $M,s\sim_n M',s'$, i.e., each $v\in s$ is $n$-bisimilar to some $v'\in s'$ and vice versa. 

In one direction, take $v\in s$. Since $M,s\models\bigvee\{\chi_{M',w'}^n\mid w'\in s'\}$, by persistency $M,v\models\bigvee\{\chi_{M',w'}^n\mid w'\in s'\}$. This implies that for some $v'\in s'$ we have $M,v\models\chi_{M',v'}^n$, which by the first claim implies
 $M,v\sim_n M',v'$. 
 
In the opposite direction, take $v'\in s'$. Since $M,s\not\models\Lori\{\neg\chi_{M',w'}^n\mid w'\in s'\}$, in particular $M,s\not\models\neg\chi_{M',v'}^n$. Since $\neg\chi_{M',v'}^n$ is a declarative, and thus truth-conditional, it follows that there is a $v\in s$ with $M,v\not\models\neg\chi_{M',v'}^n$; since truth conditions work in the usual way, $M,v\models\chi_{M',v'}^n$, which by the first claim implies $M,v\sim_n M',v'$.

\item[$\Leftarrow$] Suppose $M,s\sim_n M',s'$ for some $s\in\Sigma(w)$. We will show that $M,s\models\bigvee\{\chi_{M',w'}^n\mid w'\in s'\}$ but $M,s\not\models\Lori\{\neg\chi_{M',w'}^n\mid w'\in s'\}$, which implies $M,w\models\pi_{M',s'}^n$. 

To show that $M,s\models\bigvee\{\chi_{M',w'}^n\mid w'\in s'\}$, we reason as follows. Take any $v\in s$. Since $s\sim_n s'$, there is a $v'\in s'$ with $v\sim_n v'$. By the first claim, we have $M,v\models\chi_{M',v'}^n$, and therefore $M,v\models\bigvee\{\chi_{M',w'}^n\mid w'\in s'\}$. Since $v$ was an arbitrary world in $s$ and the formula $\bigvee\{\chi_{M',w'}^n\mid w'\in s'\}$ is truth-conditional (since it is a declarative) we have $M,s\models \bigvee\{\chi_{M',w'}^n\mid w'\in s'\}$, as desired.

To show that $M,s\not\models\Lori\{\neg\chi_{M',w'}^n\mid w'\in s'\}$, 
take an arbitrary $v'\in s'$. Since $M,s\sim_n M',s'$, there is a $v\in s$ with $M,v\sim_n M',v'$. By the first claim, $M,v\models\chi_{M',v'}^n$, and so $M,v\not\models\neg\chi_{M',v'}^n$. By persistency, also $M,s\not\models\neg\chi_{M',v'}^n$. Since this holds for an arbitrary $v'\in s'$, we have $M,s\not\models\Lori\{\neg\chi_{M',w'}^n\mid w'\in s'\}$, as desired.
\end{itemize}
Having secured the implication (1) $\Rightarrow$ (2) for each given $n$, we proceed to prove by induction that (1) holds for every $n\in\mathbb{N}$. If $n=0$ the claim is obvious. So, suppose (1) and (2) hold for $n$; we show that the biconditional in (1) holds for $n+1$.
\begin{itemize}
\item[$\Rightarrow$] Suppose $M,w\models\chi_{M',w'}^{n+1}$. By definition of $\chi_{M',w'}^{n+1}$, this means that three things hold:
\begin{itemize}
\item[(i)] $M,w\models\chi_{M',w'}^0$; 
\item[(ii)] $M,w\models\pi_{M',s'}^n$ for every $s'\in\Sigma'(w')$;
\item [(iii)] $M,w\models\neg\pi_{M'',s''}^n$ for every $\<M'',s''\>$ s.t.\ $\neg\exists s'\in\Sigma'(w')$ with $M'',s''\sim_n M',s'$. 
\end{itemize}
Note that condition (iii) is equivalent to the following:
\begin{itemize}
\item[(iii')] for every $\<M'',s''\>$, if $M,w\models\pi_{M'',s''}^n$ there is $s'\in\Sigma(w')$ with $M'',s''\sim_n M',s'$. 
\end{itemize}

Now, (i) implies that $w$ and $w'$ make true the same atoms, i.e., the atomic condition for $M,w\sim_{n+1}M',w'$ is satisfied; (ii) implies, by the induction hypothesis, that for every $s'\in\Sigma'(w')$ there is $s\in\Sigma(w)$ with $M,s\sim_n M',s'$, i.e., the back condition for $M,w\sim_{n+1}M',w'$ is satisfied. To see that the forth condition for $M,w\sim_{n+1} M',w'$ is satisfied as well, take any $s\in\Sigma(w)$. By the induction hypothesis, $M,w\models\pi_{M,s}^n$, and so by condition (iii) there is some $s'\in\Sigma'(w')$ with $M,s\sim_n M',s'$, as required.
\item[$\Leftarrow$] Suppose $M,w\sim_{n+1}M',w'$. To show that $M,w\models\chi_{M',w'}^{n+1}$, we need to establish conditions (i)-(iii) listed above. 

The atomic condition on $\sim_{n+1}$ directly implies (i). 

The back condition ensures that for all $s'\in\Sigma'(w')$ there is $s\in\Sigma(w)$ with $M,s\sim_n M',s'$, which by the induction hypothesis implies $M,w\models\pi_{M',s'}^n$, securing (ii).

Finally, to establish (iii), suppose $M'',s''$ is a state-pointed model such that there is no $s'\in\Sigma'(w')$ with $M'',s''\sim_n M',s'$. 
Now suppose towards a contradiction that there is some $s\in\Sigma(w)$ with $M'',s''\sim_n M,s$; by the forth condition on $w\sim_{n+1}w'$, there is $s'\in\Sigma'(w')$ with $M,s\sim_n M',s'$, whence by the transitivity of $\sim_n$ it follows $M'',s''\sim_n M',s'$, contrary to assumption. So, there is no $s\in\Sigma(w)$ with $M'',s''\sim_n M,s$. Therefore, the induction hypothesis implies $M,w\not\models\pi_{M'',s''}^n$, and since truth conditions work in the usual way, $M,w\models\neg\pi_{M'',s''}^n$. \hfill$\Box$
\end{itemize}

\medskip
\noindent By using the characteristic formulas $\chi_{M,w}^n$, we can show that modal equivalence on formulas of modal depth $n$ implies $n$-bisimilarity, both at the level of worlds and at the level of information states. In other words, any behavioral difference that shows up within $n$ modal transition steps can be detected by a formula of modal depth $n$.

\begin{theorem}[Over finite languages, $n$-modal equivalence implies $n$-bisimilarity]\label{th:distinguishing power} If the set of atoms \P\ is finite, we have: 
\begin{itemize}
\item $M,w\modeq_n M',w'$ implies $M,w\sim_n M',w'$
\item $M,s\modeq_n M's'$ implies $M,s\sim_n M',s'$.
\end{itemize}
\end{theorem}

\begin{proof} Suppose $M,w\modeq_n M',w'$. Since $M',w'\models\chi_{M',w'}^n$ and $\chi_{M',w'}^n$ has modal depth $n$, it follows that $M,w\models\chi_{M',w'}^n$, which by the previous lemma implies $M,w\sim_n M',w'$.

Next, consider the claim for information states. Suppose $M,s\not\sim_n M',s'$. This means that some $w\in s$ is not $n$-bisimilar to any $w'\in s'$, or vice versa. Without loss of generality, suppose the former is the case. Then by the previous lemma, for every $w'\in s'$ we have $M',w'\models\neg\chi_{M,w}^n$. Since $\neg\chi_{M,w}^n$ is a declarative, and thus truth-conditional, it follows that $M',s'\models\neg\chi_{M,w}^n$. However, $M,w\not\models\neg\chi_{M,w}^n$, and since $w\in s$, by persistence $M,s\not\models\neg\chi_{M,w}^n$. Since  $\neg\chi_{M,w}^n$ is a formula of modal depth $n$, this shows that $M,s\not\modeq_n M',s'$.
\end{proof}

\noindent
It is worth noting that the situation for worlds and states is not fully analogous: in the case of worlds, $M,w\not\sim_n M',w'$ implies that there is a formula of modal depth $n$ that is true at $w$ but not at $w'$ \emph{and} one that is true at $w'$ but not at $w$. In the case of states, $M,s\not\sim_n M',s'$ only implies that there is a formula of modal depth $n$ supported by $s$ but not $s'$, \emph{or} one supported by $s'$ but not by $s$. Indeed, even if $M,s\not\sim_n M',s'$ it may be that every formula supported by $s$ is also supported by $s'$; in particular, this happens by persistency when $M=M'$ and $s'\subset s$.

As in standard modal logic, modal equivalence does not in general guarantee full bisimilarity: by Theorem \ref{th:distinguishing power}, modal equivalence amounts to $n$-bisimilarity for each $n$, which is strictly weaker than full bisimilarity, as shown by Example \ref{ex:infinity}. However, what we can prove is an analogue of the classical Hennessy-Milner theorem. To state this analogue, we first need to adapt the notion of image-finiteness to neighborhood models.

\begin{definition}[Image-finiteness] A neighborhood model $M=\<W,\Sigma,V\>$ is image-finite if for every world $w$, the union $\bigcup\Sigma(w)$ of its neighborhoods is finite. Equivalently, $M$ is image-finite if the underlying Kripke model $M_K$ is image-finite.
\end{definition}

\noindent
Note that image-finiteness is equivalent to the conjunction of the following two conditions: (i) each world has finitely many neighborhoods and (ii) each neighborhood is finite. 

The counterpart of the Hennessy-Milner theorem can now be stated as follows.

\begin{theorem} If $M,M'$ are image-finite in-models, then: 
\begin{itemize}
\item $M,w\modeq M',w'$ implies $M,w\sim M',w'$;
\item if $s,s'$ are finite states, $M,s\modeq M',s'$ implies $M,s\sim M',s'$.
\end{itemize}
\end{theorem}

\begin{proof} Suppose $M$ and $M'$ are image-finite. We will show that the modal equivalence relation between worlds in these models, ${\modeq}\subseteq W\times W'$, is a bisimulation. 

Suppose $w\modeq w'$. Clearly, $w$ and $w'$ satisfy the same atoms. We show that the Forth condition is satisfied, and leave the analogous proof of the Back condition to the reader. 

Take a state $s\in\Sigma(w)$. We need to show that $s\modeq s'$ for some $s'\in\Sigma'(w')$. 
By image-finiteness we can write $s=\{w_1,\dots,w_n\}$ and $\bigcup\Sigma'(w')=\{v_1,\dots,v_m\}$. Note that $n\ge 1$ since we require neighborhoods to be nonempty, and $m\ge 1$ since if $m=0$ we would have $\Sigma'(w')=\emptyset$, in which case $w$ and $w'$ would be distinguished by the formula $(\top\To\bot)$, contradicting $w\modeq w'$.

For each $i\le n,j\le m$ we define a formula $\delta_{ij}$ as follows: if $w_i\modeq v_j$ we let $\delta_{ij}=\top$, while if $w_i\not\modeq v_j$ we let $\delta_{ij}$ be a declarative such that $M,w_i\models\delta_{ij}$ and $M,v_j\models\neg\delta_{ij}$.\footnote{Note that such a declarative exists: if $w_i\not\modeq v_j$, then $w_i$ and $v_j$ disagree about the truth of some formula $\xi\in\L$. Now $\xi$ has the same truth conditions as the declarative $\xi^!$ defined as in the proof of Proposition \ref{prop:declarative variant}, so $w_i$ and $v_i$ disagree on the truth of $\xi^!$. We can choose $\delta=\xi^!$ if $w_i\models\xi^!$ and $\delta=\neg\xi^!$ otherwise.}
If we now define $\gamma_i:=\bigwedge_{j\le m}\delta_{ij}$, we have that for all $j\le m$: 
$$(*)\qquad M',v_j\models\gamma_{i}\iff w_i\modeq v_j\qquad\phantom{(*)}$$
To see this, note that all conjuncts of $\gamma_i$ are true at $w_i$ and so if $w_i\modeq v_j$, they are true at $v_j$ as well; on the other hand, if $w_i\not\modeq v_j$, then by construction $\gamma_i$ includes a conjunct $\delta_{ij}$ which is false at $v_j$. 

Now consider the formula:
$$\aa\;:=\;(\bigvee_{i=1}^{n}\gamma_i)\To(\Lori_{i=1}^n\neg\gamma_i)$$
It is easy to see that $s$ supports the antecedent of~$\aa$ but not the consequent, so $M,w\not\models\aa$. 
Since $w\modeq w'$, we have that $M',w'\not\models\aa$. So there is a $s'\in\Sigma'(w')$ that supports the antecedent of $\aa$ but not the consequent. We are now going to show that $s\modeq s'$, i.e., every world in $s$ is modally equivalent to a world in $s'$ and vice versa.

Since $s'\models\bigvee_{i\le n}\gamma_i$, by persistency every world $v\in s'$ satisfies some formula $\gamma_i$ and thus by $(*)$ it is modally equivalent to some world in $s$, namely $w_i$. For the converse, take a world $w_i\in s$. Since $s'$ does not support the consequent of $\aa$ we have $s'\not\models\neg\gamma_i$. Since $\neg\gamma_i$ is a declarative and thus truth-conditional, there is a world $v\in s'$ such that $v\not\models\neg\gamma_i$ and so $v\models\gamma_i$, which by $(*)$ implies $w_i\modeq v$. 
This proves the claim for worlds. 

For the claim about states, let $s,s'$ be finite states with $s\modeq s'$. Since $\modeq$ on worlds is a bisimulation, it suffices to show that every world in $s$ is modally equivalent to a world in $s'$ and vice~versa. So, take $w\in s$. 
Proceeding in the same way as above, we can define a declarative $\gamma_w$ such that $w\models\gamma_w$ and for all $v\in s'$ we have $v\models\gamma_w\iff w\modeq v$. Now since $w\in s$, by persistency we have $s\not\models\neg\gamma_w$ and since $s\modeq s'$ also $s'\not\models\neg\gamma_w$. Since $\neg\gamma_w$ is a declarative and thus truth-conditional, there is some $w'\in s'$ with $w'\not\models\neg\gamma_w$ and so $w'\models\gamma_w$, which implies $w\modeq w'$. The converse is proved analogously.
\end{proof}

\noindent Note that the condition that $s,s'$ be \emph{finite} states is crucial: it is easy to construct a counterexample in the spirit of Example \ref{ex:infinity}, involving an image-finite model $M$ and two infinite states $\{w_0,w_1,w_2,\dots\}$ and $\{w_0,w_1,w_2,\dots, w_{\omega}\}$ which are $n$-bisimilar for each $n$ (since $M,w_\omega\sim_n M,w_n$) but not fully bisimilar (since $M,w_\omega\not\sim M,w_n$ for any $n$).

We conclude this section with a characterization of the properties of worlds and states 
that can be defined in \inqnl. 

\begin{definition}[World and state properties] By a \emph{world-property} (or \emph{state-property}) we mean a class $\C$ of world-pointed (or state-pointed) in-models. 

We say that a world-property $\C$ is:
\begin{itemize}
\item \emph{closed under $\sim_n$}, if $\<M,w\>\in\C$ and  $M,w\sim_n M',w'$ implies $\<M',w'\>\in\C$;
\item \emph{definable in \inqnl} if for some formula $\phi\in\L_\To$, $\C=\{\<M,w\>\mid M,w\models\phi\}$.
\end{itemize}
Analogous definitions apply to state-properties. Furthermore, if \C\ is a state-property, we say that \C\ is 
\emph{downward-closed}, if $\<M,s\>\in\C$ and  $t\subseteq s$ implies $\<M,t\>\in\C$.
\end{definition}

\begin{theorem}[Expressive power of \inqnl\ with respect to worlds]\label{th:expressive power worlds}
Provided the set of atoms is finite, a world-property is definable in \inqnl\ iff it is $\sim_n$-closed for some $n$.
\end{theorem}

\begin{proof} The left-to-right direction follows from the fact that any formula $\phi$ of \inqnl\ is preserved under $\sim_n$ for any $n\ge \md(\phi)$. For the right-to-left direction, consider a world-property \C\ and suppose \C\ is closed under $\sim_n$. Consider the formula:
$$\rho_\C:=\bigvee\{\chi_{M,w}^n\mid \<M,w\>\in\C\}$$
This formula is well-defined, since we argued in the proof of Lemma \ref{lemma:characteristic formulas} that for a fixed $n$ there are only finitely many distinct formulas of the form $\chi_{M,w}^n$. We claim that \C\ is defined by $\rho_\C$, i.e., for any world-pointed model $\<M,w\>$ we have
$$\<M,w\>\in\C\iff M,w\models\rho_\C$$ 
For the $\Rightarrow$ direction, if $\<M,w\>\in\C$ then $\chi_{M,w}^n$ is a disjunct of $\rho_\C$. By Lemma \ref{lemma:characteristic formulas} we have $M,w\models\chi_{M,w}^n$ and so also $M,w\models\rho_\C$. 

For the converse, suppose $M,w\models\rho_\C$. Then $M,w\models\chi_{M',w'}^n$ for some $\<M',w'\>\in\C$. By Lemma \ref{lemma:characteristic formulas}, $M,w\sim_n M',w'$, and since $\C$ is closed under $\sim_n$, also $\<M,w\>\in\C$.
\end{proof}

\begin{theorem}[Expressive power of \inqnl\ with respect to information states]\label{th:expressive power states} Provided the set of atoms is finite, a state-property \C\ is definable in \inqnl\ iff \C\ is $\sim_n$-closed for some $n$, non-empty, and downward closed.
\end{theorem}

\begin{proof} The left-to-right direction follows from Proposition \ref{prop:n-bis to n-modeq} together with the persistence of support and the empty state property. For the converse, let $\C\neq\emptyset$ be a state-property which is $\sim_n$-invariant for a given $n$ and downward closed. We are going to show that \C\ is definable in \inqnl.

For a state-pointed model $\<M,s\>$, define $\chi^n_{M,s}=\bigvee\{\chi^n_{M,w}\mid w\in s\}$ (if $s=\emptyset$, we let $\chi^n_{M,s}=\bot$). It is straightforward to verify that we have:
$$M,s\models \chi^n_{M',s'}\iff M,s\sim_n M',t'\text{ for some }t'\subseteq s'$$
Now consider the following formula:
$$\sigma_\C:=\Lori\{\chi_{M,s}^n\mid \<M,s\>\in\C\}$$
We claim that $M,s\models\sigma_\C\iff \<M,s\>\in \C$. For the right-to-left direction, if $\<M,s\>\in\C$, then $\chi^n_{M,s}$ is a disjunct of $\sigma_\C$; since $\<M,s\>$ supports $\chi^n_{M,s}$, it supports $\sigma_\C$. For the converse, suppose $M,s\models\sigma_\C$. This means that for some $\<M',s'\>\in\C$ we have $M,s\models\chi^n_{M',s'}$, which implies that $M,s\sim_n M',t'$ for some $t'\subseteq s'$. Since \C\ is downward-closed we have $\<M',t'\>\in\C$, and finally by $\sim_n$-closure, $\<M,s\>\in\C$.
\end{proof}

\section{Axiomatization}
\label{sec:axiomatization}

\noindent Having explored the expressive power of \inqnl, we now turn to its axiomatization. In this section, we provide a sound and complete Hilbert-style proof system for the logic \inqnl, as interpreted over the class of all in-models. In the next section we then turn to the logic of some specific classes of models. 
The axioms of our proof system can be divided into two classes: propositional axioms, capturing the logical properties of the connectives and the logical relation between statements and questions, and modal axioms, capturing the logical properties of $\To$. The propositional axioms are: 

\begin{itemize}
\item all instances of the axioms for intuitionistic propositional logic, with $\lori$ in the role of intuitionistic disjunction;
\item all instances of the following two schemata, where $\aa$ ranges only over declaratives, while $\phi$ and $\psi$ range over arbitrary formulas:

\medskip\noindent
\begin{tabular}{ll}
(\texttt{DDN}) & $\neg\neg\aa\to\aa$\\
(\texttt{Split}) & $(\aa\to\phi\lori\psi)\to(\aa\to\phi)\lori(\aa\to\psi)$
\end{tabular}

\end{itemize}

\medskip\noindent
Intuitively, (\texttt{DDN}) captures the fact that declaratives obey classical logic, while (\texttt{Split}) captures the fact that they are statements \cite[cf.][for detailed discussion]{Ciardelli:23book}.

The modal axioms are all instances of the following three schemata:

\medskip\noindent
\begin{tabular}{ll}
(\texttt{Trans}) & $(\phi\To\psi)\land(\psi\To\chi)\to(\phi\To\chi)$\\
(\texttt{RConj}) & $(\phi\To\psi)\land(\phi\To\chi)\to(\phi\To \psi\land\chi)$\\
(\texttt{LDisj}) & $(\phi\To\chi)\land(\psi\To\chi)\to(\phi\lori\psi\To\chi)$
\end{tabular}

\medskip\noindent The inference rules of the system are modus ponens and a conditional version of necessitation:
$$\infer[(\texttt{MP})]{\psi}{\phi & \phi\to\psi}\qquad\qquad \infer[(\texttt{CN})]{\phi\To\psi}{\phi\to\psi}$$
\noindent If $\Phi,\Psi\subseteq\L$, we write $\Phi\vdash\Psi$ if there are $\phi_1,\dots,\phi_n\in\Phi$ and $\psi_1,\dots,\psi_m\in\Psi$ such that the formula $(\phi_1\land\dots\land\phi_n)\to(\psi_1\lori\dots\lori\psi_m)$ is derivable in the system (we allow for $n=0$, in which case the relevant antecedent is $\top$, and for $m=0$, in which case the consequent is $\bot$). We write $\phi_1,\dots,\phi_n\vdash\psi_1,\dots,\psi_m$ instead of $\{\phi_1,\dots,\phi_n\}\vdash\{\psi_1,\dots,\psi_m\}$, and we write $\phi\vdasheq\psi$ in case $\phi\vdash\psi$ and $\psi\vdash\phi$.

As usual, the soundness of the system is established by checking that the axioms are valid and the rules preserve validity, which is straightforward.

\begin{proposition}[Soundness] For all $\Phi\cup\{\psi\}\subseteq\L_\To$, $\Phi\vdash\psi$ implies $\Phi\models\psi$.
\end{proposition}

\noindent Towards the proof of completeness, 
we start with a few preliminary results concerning resolutions, which are not specific to \inqnl, but hold for inquisitive propositional logic and its modal extensions in general. 

First, our proof system suffices to prove the equivalence between a formula $\phi$ and its normal form given by Proposition \ref{prop:normal form}.

\begin{lemma}\label{lemma:provable NF} For all $\phi\in\L_\To$, $\phi\vdasheq\Lori\R(\phi)$.
\end{lemma}

\noindent
This fact depends only on the propositional component of the proof system, and can be proved in exactly the same way as for inquisitive propositional logic: see Lemma 4.3.4 in \cite{Ciardelli:23book} for the details. 
As an immediate corollary we have the following fact.

\begin{lemma}\label{lemma:resolutions prove formula} For all $\phi\in\L_\To$ and all $\aa\in\R(\phi)$, $\aa\vdash\phi$.
\end{lemma}

\noindent The next lemma says that if a set of declaratives derives a formula, it derives a particular resolution of it (note that the converse holds by~the previous lemma). The proof is again standard in inquisitive logic \citep[cf.\ Lemma 4.3.8 in][]{Ciardelli:23book}. 

\begin{lemma}\label{lemma:provable split} Let $\Gamma\subseteq\L_\To^!$ and $\phi\in\L_\To$. If $\Gamma\vdash\phi$, then $\Gamma\vdash\aa$ for some $\aa\in\R(\phi)$.
\end{lemma}

\noindent For our result, we will also need a notion of resolutions for sets of formulas. A resolution of a set $\Phi$ is a set obtained by replacing each element $\phi\in\Phi$ by a resolution. Thus, for example, the resolutions of the set $\Phi=\{p,?q\}$ are $\Gamma_1=\{p,q\}$ and $\Gamma_2=\{p,\neg q\}$. 

\begin{definition} A resolution function for a set of formulas $\Phi\subseteq\L_\To$ is a function $f$ that associates to each $\phi\in\Phi$ some resolution $f(\phi)\in\R(\phi)$. A resolution of $\Phi$ is  the image of $\Phi$ under a resolution function: $\R(\Phi)=\{f[\Phi]\mid f\text{ a resolution function for }\Phi\}$.
\end{definition}

\noindent Note that if $\Gamma$ is a set of declaratives then $\R(\Gamma)=\{\Gamma\}$.

The next lemma says that if $\Phi$ fails to derive $\Psi$, then $\Phi$ can be strengthened to a resolution $\Gamma\in\R(\Phi)$ that still fails to derive $\Psi$. The proof is again the same as for inquisitive propositional logic: see Lemma 4.3.7 in \cite{Ciardelli:23book} for the details.

\begin{lemma}\label{lemma:provable specification} For all $\Phi,\Psi\subseteq\L$, if $\Phi\not\vdash\Psi$ then there is $\Delta\in\R(\Phi)$ s.t.\ $\Delta\not\vdash\Psi$.
\end{lemma}

\noindent
With these standard preliminary results at hand, we are now ready to describe a canonical model construction for $\inqnl$. This construction is based on \emph{complete theories of declaratives} (CTDs), defined as follows.

\begin{definition} A set $\Gamma\subseteq\L_\To^!$ of declaratives is a \emph{complete theory of declaratives} (CTD for short) if it satisfies the following conditions:
\begin{enumerate}
\item  $\Gamma$ is deductively closed w.r.t.\ declaratives: if $\Gamma\vdash\aa$ and $\aa\in\L_\To^!$ then $\aa\in\Gamma$;
\item $\Gamma$ is consistent: $\bot\not\in\Gamma$;
 \item $\Gamma$ is complete: for all $\aa\in\L_\To^!$, one of $\aa$ and $\neg\aa$ is in $\Gamma$. 
\end{enumerate}
\end{definition}
\noindent
If $S$ is a set of CTDs, we let $\bigcap S=\{\aa\in\L_\To^!\mid \aa\in\Gamma\text{ for all }\Gamma\in S\}$ (thus,  $\bigcap\emptyset=\L_\To^!$). Note that $\bigcap S$ is deductively closed with respect to declaratives, i.e., the following holds.

\begin{remark} \label{remark:bigcapS}
For any set $S$ of CTDs and any $\aa\in\L_\To^!$, $\bigcap S\vdash\aa$ implies $\aa\in\bigcap S$. 
\end{remark}

\begin{proof} Note that $\bigcap S\subseteq\Gamma$ for any $\Gamma\in S$. Using this fact and the fact that CTDs are closed under deduction of declaratives, we have:
\begin{eqnarray*}
\bigcap S\vdash\aa&\Rightarrow& \text{for all }\Gamma\in S: \Gamma\vdash\aa\\
&\Rightarrow& \text{for all }\Gamma\in S: \aa\in\Gamma\\
&\Rightarrow& \aa\in\bigcap S
\end{eqnarray*}
 
\vspace{-.75cm}
\end{proof}

\noindent
With any set $\Delta\subseteq\L_\To^!$ of declaratives we can associate a set of CTDs, namely, the set of its complete extensions: $S_\Delta=\{\Gamma\mid\Gamma\text{ is a CTD and }\Delta\subseteq\Gamma\}$. The following Lindenbaum-type lemma, which says that a consistent theory of declaratives has a complete extension, is proved by the standard saturation argument.

\begin{lemma}\label{lemma:extension} If $\Delta\subseteq\L_\To^!$ and $\Delta\not\vdash\bot$, then $S_\Delta\neq\emptyset$.
\end{lemma}

\noindent For any set of declaratives $\Delta\subseteq\L_\To^!$, the sets $\Delta$ and $\bigcap S_\Delta$ prove the same formulas.

\begin{lemma}\label{lemma:intersection} For any $\Delta\subseteq\L_\To^!$ and $\phi\in\L$: $\Delta\vdash\phi\iff\bigcap S_\Delta\vdash\phi$.
\end{lemma}

\begin{proof} The direction $\Rightarrow$ is clear as $\Delta\subseteq\bigcap S_\Delta$. For the converse, suppose for a contradiction that for some $\phi$ we had $\bigcap S_\Delta\vdash\phi$ but $\Delta\not\vdash\phi$. Since $\bigcap S_\Delta$ is a set of declaratives, by Lemma \ref{lemma:provable split} we have $\bigcap S_\Delta\vdash \aa$ for some $\aa\in\R(\phi)$. Since $\Delta\not\vdash\phi$, it follows by Lemma \ref{lemma:resolutions prove formula} that $\Delta\not\vdash\aa$. By the axiom $\neg\neg\aa\to\aa$, this implies $\Delta\not\vdash\neg\neg\aa$, and therefore $\Delta,\neg\aa\not\vdash\bot$.
So by Lemma \ref{lemma:extension} there is a CTD $\Gamma'$ such that $\Delta\cup\{\neg\aa\}\subseteq\Gamma'$. Since $\Gamma'\in S_\Delta$ and $\aa\not\in\Gamma'$ we have $\aa\not\in\bigcap S_\Delta$, so by Remark \ref{remark:bigcapS} we have $\bigcap S_\Delta\not\vdash\aa$, which is a contradiction.
\end{proof}

\noindent
We now define our canonical model as follows. 

\begin{definition} The canonical model for \inqnl\ is $M^c=\<W^c,\Sigma^c,V^c\>$ where:
\begin{itemize}
\item $W^c$ is the set of complete theories of declaratives;
\item ${\Sigma^c(\Gamma)=\{S\neq\emptyset\mid \forall \phi,\psi\in\L_\To:(\phi\To\psi)\in\Gamma\text{ and }\bigcap S\vdash\phi\text{ implies }\bigcap S\vdash\psi\}}$
\item $V^c(p)=\{\Gamma\in W^c\mid p\in\Gamma\}$
\end{itemize}
\end{definition}

\noindent
The key to the completeness proof is an analogue of the standard existence lemma. 
\begin{lemma}[Existence Lemma]\label{lemma:existence} If $\Gamma$ is a CTD and $(\phi\To\psi)\not\in\Gamma$, there exists a state $S\in\Sigma^c(\Gamma)$ such that $\bigcap S\vdash\phi$ and $\bigcap S\not\vdash\psi$.  
\end{lemma}

\noindent
Towards the proof of this lemma, we first establish some preliminary results. Let $\Gamma$ be a fixed CTD. Given two sets $\Phi,\Psi\subseteq\L_\To$, we write ${\Phi\Togamma\Psi}$ if there are $\phi_1,\dots,\phi_n\in\Phi$ and $\psi_1,\dots,\psi_m\in\Psi$ such that the formula $((\bigwedge_{i\le n}\phi_i)\To(\Lori_{i\le m}\psi_i))$ is in $\Gamma$. The following cut-like property of the relation $\Togamma$ is crucial for our construction.

\begin{lemma}\label{lemma:cut} For any two sets $\Phi,\Psi\subseteq\L_\To$ and any formula $\chi\in\L_\To$: 
$\Phi\cup\{\chi\}\Togamma\Psi$ and $\Phi\Togamma\Psi\cup\{\chi\}$ implies $\Phi\Togamma\Psi$.
\end{lemma}

\begin{proof}
Suppose $\Phi\cup\{\chi\}\Togamma\Psi$ and $\Phi\Togamma\Psi\cup\{\chi\}$. This means that there are formulas ${\phi_1,\dots,\phi_n,\phi_{n+1},\dots,\phi_{n+m}}\in\Phi$ and $\psi_1,\dots,\psi_h,\psi_{h+1},\dots,\psi_{h+k}\in\Psi$ such that:
$$((\chi\land\bigwedge_{1=i}^n\phi_i)\To\Lori_{1=i}^h\psi_i)\in\Gamma\qquad((\bigwedge_{i=1}^{m}\phi_{n+i})\To(\chi\lori\Lori_{i=1}^{k}\psi_{h+i}))\in\Gamma$$
We show that $\Gamma$ contains $((\bigwedge_{i\le n+m}\phi_i)\To(\Lori_{i\le h+k}\psi_i))$, thus witnessing ${\Phi\Togamma\Psi}$. To ease notation, we spell out the details for the case $n=m=h=k=1$, but the general case is completely analogous. It suffices to show that:
$$(\phi_1\land\chi\To\psi_1),(\phi_2\To\psi_2\lori\chi)\vdash(\phi_1\land\phi_2\To\psi_1\lori\psi_2)$$
Since the formulas on the left-hand-side are in $\Gamma$ and $\Gamma$ is closed under deduction of declaratives, so is the conclusion. 

First note that the following formula is provable in the propositional component of the proof system using the standard axioms for $\land$ and $\lori$:
\begin{equation}\label{distribution}
\phi_1\land(\psi_2\lori\chi)\to\psi_2\lori(\phi_1\land\chi)
\end{equation}
In the following derivation, we indicate explicitly only the modal axioms and rules involved in the reasoning, omitting reference to propositional axioms. 
For simplicity, we use the formulas $(\phi_1\land\chi\To\psi_1)$ and $(\phi_2\To\psi_2\lori\chi)$ as if they were premises; this is legitimate since we will not use the conditional necessitation rule (\texttt{CN}) on these formulas or anything inferred from them. Rewriting the argument with the relevant formulas used throughout as conditional antecedents is tedious but straightforward. 

\medskip\noindent
\begin{tabular}{lll}
1. & $\phi_1\land\chi\,\To\,\psi_1$ & (premise)\\
2. & $\phi_2\,\To\,\psi_2\lori\chi$ & (premise)\\
3. & $\phi_1\land\phi_2\To\phi_1$ & (\texttt{CN}) from axiom $\phi_1\land\phi_2\to\phi_1$\\
4. & $\phi_1\land\phi_2\To\phi_2$ & (\texttt{CN})  from axiom $\phi_1\land\phi_2\to\phi_2$\\
5. & $\phi_1\land\phi_2\To\psi_2\lori\chi$ & (\texttt{MP}) from (\texttt{Trans}), 4, 2\\
6. & $\phi_1\land\phi_2\To\phi_1\land(\psi_2\lori\chi)$ & (\texttt{MP}) from (\texttt{RConj}), 3, 5\\
7. & $\phi_1\land(\psi_2\lori\chi)\To\psi_2\lori (\phi_1\land\chi)$ & (\texttt{CN}) from (\ref{distribution})\\
8. & $\phi_1\land\phi_2\To\psi_2\lori (\phi_1\land\chi)$ & (\texttt{MP}) from (\texttt{Trans}), 6, 7\\
9. & $\psi_1\To\psi_1\lori\psi_2$ & (\texttt{CN}) from axiom $\psi_1\to\psi_1\lori\psi_2$\\
10. & $\psi_2\To\psi_1\lori\psi_2$ & (\texttt{CN}) from axiom $\psi_2\to\psi_1\lori\psi_2$\\
11. & $\phi_1\land\chi\To\psi_1\lori\psi_2$ & (\texttt{MP}) from (\texttt{Trans}), 1, 9\\
12. & $\psi_2\lori (\phi_1\land\chi)\To \psi_1\lori\psi_2$ & (\texttt{MP}) from (\texttt{LDisj}), 10, 11\\
13. & $\phi_1\land\phi_2\,\To\,\psi_1\lori\psi_2$ & (\texttt{MP}) from (\texttt{Trans}), 8, 12
\end{tabular}

\vspace{-.4cm}
\end{proof}

\smallskip
\begin{lemma}[Splitting lemma]\label{lemma:splitting} Let $\Gamma$ be a CTD and suppose $\Phi\not\Togamma\Psi$. Then we can partition the language $\L_\To$ into two sets $\Left$ and $\Right$ such that $\Phi\subseteq\Left$, $\Psi\subseteq\Right$, and $\Left\not\Togamma\Right$.
\end{lemma}

\begin{proof}
Fix an enumeration $(\chi_n)_{n\in\mathbb N}$ of $\L_\To$. 
We define a sequence of sets $(\Left_n)_{n\in\mathbb{N}}$ and $(\Right_n)_{n\in\mathbb{N}}$ as follows:
\begin{itemize}
\item $\Left_0=\Phi, \Right_0=\Psi$
\item if $\Left_n\cup\{\chi_n\}\not\Togamma\Right_n$ we let $\Left_{n+1}:=\Left_n\cup\{\chi_n\}$ and $\Right_{n+1}=\Right_n$
\item if $\Left_n\cup\{\chi_n\}\Togamma\Right_n$ we let $\Left_{n+1}:=\Left_n$ and $\Right_{n+1}=\Right_n\cup\{\chi_n\}$
\end{itemize}
We show by induction on $n$ that $\Left_n\not\Togamma\Right_n$. For $n=0$ this is true by assumption. Now suppose this is true for $n$ and consider $n+1$. If $\Left_n\cup\{\chi_n\}\not\Togamma\Right_n$, the claim is obvious from the definition. So, suppose $\Left_n\cup\{\chi_n\}\Togamma\Right_n$. Since by induction hypothesis $\Left_n\not\Togamma\Right_n$, Lemma \ref{lemma:cut} implies $\Left_n\not\Togamma\Right_n\cup\{\chi_n\}$, which by definition amounts to $\Left_{n+1}\not\Togamma\Right_{n+1}$.

Now let $\Left=\bigcup_{n\in\mathbb{N}}\Left_n$ and $\Right=\bigcup_{n\in\mathbb{N}}\Right_n$. By construction, $\Phi\subseteq\Left$ and $\Psi\subseteq\Right$. We have $\Left\not\Togamma\Right$, since otherwise there would be an $n$ such that $\Left_n\Togamma\Right_n$. Moreover, $\Left$ and $\Right$ form a partition of $\L$: by construction, every formula occurs in either set, and no formula can occur in both (if $\chi\in\Left\cap\Right$, then since $(\chi\To\chi)\in\Gamma$ we would have $\Left\Togamma\Right$).
\end{proof}

\noindent
We are now finally ready for the proof of the existence lemma.

\medskip\noindent\emph{Proof of Lemma \ref{lemma:existence}.}
Let $\Gamma$ be a CTD with $(\phi\To\psi)\not\in\Gamma$. This means that $\{\phi\}\not\Togamma\{\psi\}$.
Extend $\{\phi\}$ and $\{\psi\}$ to sets $\Left$ and $\Right$ as in the previous lemma. 
Since $\Left\not\Togamma\Right$, the rule (\texttt{CN}) guarantees that $\Left\not\vdash\Right$. 
By Lemma \ref{lemma:provable specification} we can find a set $\Delta\in\R(\Left)$ with $\Delta\not\vdash\Right$. We can now take $S=S_\Delta=\{\Gamma'\in W^c\mid \Delta\subseteq\Gamma'\}$. 
We need to verify that (i) $\bigcap S\vdash\phi$, (ii) $\bigcap S\not\vdash\psi$ and (iii) $S\in\Sigma^c(\Gamma)$.

\begin{itemize}
\item For (i), recall that $\phi\in\Left$. Since $\Delta\in\R(\Left)$, for some $\aa\in\R(\phi)$ we have $\aa\in\Delta$. By Lemma \ref{lemma:resolutions prove formula}, $\Delta\vdash\phi$, and thus by 
Lemma \ref{lemma:intersection} also $\bigcap S\vdash\phi$.

\item For (ii), recall that $\psi\in\Right$. Since $\Delta\not\vdash\Right$ we have $\Delta\not\vdash\psi$. By Lemma \ref{lemma:intersection}, $\bigcap S\not\vdash\psi$.

\item For (iii), first note that since $\Delta\not\vdash\Right$, we have $\Delta\not\vdash\bot$, so by Lemma \ref{lemma:extension}, $S\neq\emptyset$. Next, suppose $(\chi\To\xi)\in\Gamma$ and $\bigcap S\vdash\chi$. We need to show that $\bigcap S\vdash\xi$. By Lemma \ref{lemma:intersection}, $\Delta\vdash\chi$. Since by construction $\Delta\not\vdash\Right$, it follows that $\chi\not\in\Right$, and since $\Right$ and $\Left$ partition the set of formulas, $\chi\in\Left$. Now we must have $\xi\in\Left$ as well, for if we had $\xi\in\Right$ it would follow from $(\chi\To\xi)\in\Gamma$ that $\Left\Togamma\Right$, contrary to the properties of these sets. Since $\xi\in\Left$ and $\Delta\in\R(\Left)$, for some $\aa\in\R(\xi)$ we have $\aa\in\Delta$, so by Lemma \ref{lemma:resolutions prove formula}, $\Delta\vdash\xi$. Finally, Lemma \ref{lemma:intersection} gives $\bigcap S\vdash\xi$, as desired.\hfill$\Box$
\end{itemize}

\noindent It is interesting to remark that the Splitting Lemma also allows us to characterize the accessibility relation $R_{\Sigma^c}$ induced by the canonical model, defined by letting $R_{\Sigma^c}[\Gamma]=\bigcup\Sigma^c(\Gamma)$. Though this characterization is not used in the completeness proof itself, it will be useful in the next section.

\begin{lemma}\label{lemma:accessibility} For any CTDs $\Gamma$ and $\Gamma'$:
$$\Gamma R_{\Sigma^c}\Gamma'\iff \forall\aa\in\L_\To^!: \ibox\aa\in\Gamma\text{ implies }\aa\in\Gamma'$$
\end{lemma}

\begin{proof} Suppose $\Gamma R_{\Sigma^c}\Gamma'$, i.e., $\Gamma'\in\bigcup\Sigma^c(\Gamma)$. Then $\Gamma'\in S$ for some $S\in\Sigma^c(\Gamma)$. Now let $\aa$ be a declarative and suppose $\ibox\aa\in\Gamma$, that is, $(\top\To\aa)\in\Gamma$. Since $S\in\Sigma^c(\Gamma)$ and $\bigcap S\vdash\top$, it follows that $\bigcap S\vdash\aa$, whence by Remark \ref{remark:bigcapS}, $\aa\in\bigcap S$. Since $\Gamma'\in S$, it follows that $\aa\in\Gamma'$.

Conversely, suppose that for all $\aa\in\L_\To^!$, $\ibox\aa\in\Gamma\text{ implies }\aa\in\Gamma'$. We must show that $\Gamma'\in S$ for some $S\in\Sigma^c(\Gamma)$.

First, we claim that $\emptyset\not\Togamma\{\neg\aa\mid\aa\in\Gamma'\}$. For suppose not: then there are $\aa_1,\dots,\aa_n\in\Gamma'$ such that $(\top\To\neg\aa_1\lori\dots\lori\neg\aa_n)\in\Gamma$. Since $\neg\aa_1\lori\dots\lori\neg\aa_n\vdash\neg(\aa_1\land\dots\land\aa_n)$, by $(\texttt{CN})$ and $(\texttt{Trans})$ we have $(\top\To\neg(\aa_1\land\dots\land\aa_n))\in\Gamma$, that is, $\ibox\neg(\aa_1\land\dots\land\aa_n)\in\Gamma$. It follows that $\neg(\aa_1\land\dots\land\aa_n)\in\Gamma'$, but this is impossible since each $\aa_i$ is in $\Gamma'$ and $\Gamma'$ is consistent.

So, we have $\emptyset\not\Togamma\{\neg\aa\mid\aa\in\Gamma'\}$. By Lemma \ref{lemma:splitting}, we can partition $\L_\To$ into sets $\Left,\Right$ with $\Left\not\Togamma\Right$ and $\{\neg\aa\mid\aa\in\Gamma'\}\subseteq\Right$. Reasoning as in the previous lemma, we can find a 
$\Delta\in\R(\Left)$ with $\Delta\not\vdash\Right$, and we can show that the corresponding set of complete extensions $S_\Delta$ is in $\Sigma^c(\Gamma)$. We now claim that $\Gamma'\in S_\Delta$. To show this, it suffices to show that $\Delta\cup\Gamma'\not\vdash\bot$: if this holds, it follows by Lemma \ref{lemma:extension} that there is a CTD $\Gamma''$ with $\Delta\cup\Gamma'\subseteq\Gamma''$; since CTDs cannot be properly included in one another, we must have $\Gamma'=\Gamma''$, and therefore $\Delta\subseteq\Gamma'$, showing that $\Gamma'\in S_\Delta$ as desired.

So, towards a contradiction, suppose $\Delta\cup\Gamma'\vdash\bot$. Since $\Gamma'$ is closed under conjunction, this means that there is a formula $\aa\in\Gamma'$ such that $\Delta\cup\{\aa\}\vdash\bot$, and so, $\Delta\vdash\neg\aa$. But this is impossible, since by construction $\neg\aa\in\Right$ and $\Delta\not\vdash\Right$. 

To conclude, we have found a state $S_\Delta$ such that $\Gamma'\in S_\Delta$ and $S_\Delta\in\Sigma^c(\Gamma)$, showing that $\Gamma R_{\Sigma^c}\Gamma'$ as required.
\end{proof}

The bridge between derivability in our proof system and semantics in $M^c$ is given by the following support lemma, which generalizes the usual truth lemma.

\begin{lemma}[Support Lemma]For all states $S\subseteq W^c$ and all $\phi\in\L_\To$:\\ $M^c,S\models\phi\iff\bigcap S\vdash\phi$.
\end{lemma}

\begin{proof} By induction on $\phi$. The cases for atoms and connectives are standard (cf.\ pp.\ 90-91 in \cite{Ciardelli:16}). We spell out the inductive step for $\phi=(\psi\To\chi)$. 

Suppose $\bigcap S\vdash(\psi\To\chi)$. Take a world $\Gamma\in S$ and a state $T\in\Sigma^c(\Gamma)$ with $M^c,T\models\psi$. By induction hypothesis we have $\bigcap T\vdash\psi$. Since $\Gamma\in S$ we have $\bigcap S\subseteq\Gamma$, so $\Gamma\vdash(\psi\To\chi)$. Since $(\psi\To\chi)$ is a declarative, it follows that $(\psi\To\chi)\in\Gamma$. By definition of $\Sigma^c$, from $(\psi\To\chi)\in\Gamma$ and $\bigcap T\vdash\psi$ we can conclude $\bigcap T\vdash\chi$, which by induction hypothesis gives $M^c,T\models\chi$. Hence, $M^c,S\models(\psi\To\chi)$.

For the converse, suppose $\bigcap S\not\vdash(\psi\To\chi)$. Then there is some $\Gamma\in S$ such that $(\psi\To\chi)\not\in\Gamma$. By the Existence Lemma (Lemma \ref{lemma:existence}) there is a state $T\in\Sigma^c(\Gamma)$ such that $\bigcap T\vdash\psi$ and $\bigcap T\not\vdash\chi$, which by induction hypothesis means that $M^c,T\models\psi$ and $M^c,T\not\models\chi$. Hence, $M^c,S\not\models(\psi\To\chi)$.
\end{proof}

\noindent
Finally, we use this lemma to establish the strong completeness of our proof system.

\begin{theorem}[Strong completeness] For all $\Phi\cup\{\psi\}\subseteq\L_\To$, $\Phi\models\psi$ implies $\Phi\vdash\psi$.
\end{theorem}

\begin{proof} Suppose $\Phi\not\vdash\psi$. By Lemma \ref{lemma:provable specification} we can find a $\Delta\in\R(\Phi)$ with $\Delta\not\vdash\psi$. Note that since $\Delta\in\R(\Phi)$, for all $\phi\in\Phi$  there is some $\aa\in\R(\phi)$ with $\aa\in\Delta$, which by Lemma \ref{lemma:resolutions prove formula} gives $\Delta\vdash\phi$.
Now take $S_\Delta=\{\Gamma'\in W^c\mid \Delta\subseteq\Gamma'\}$. By Lemma \ref{lemma:intersection}, $\bigcap S_\Delta\vdash\phi$ for all $\phi\in\Phi$, while $\bigcap S_\Delta\not\vdash\psi$. By the support lemma, in the model $M^c$ the state $S_\Delta$ supports all formulas in $\Phi$ but not $\psi$, so $\Phi\not\models\psi$.
\end{proof}

\section{Frame conditions}
\label{sec:frame conditions}

\begin{figure}[t]
\begin{center}
\begin{tabular}{llll}
Condition & Definition & Axiom\\
\hline
$\downarrow$-monotonicity & $\emptyset\neq t\subseteq s\in \Sigma(w)\;\Rightarrow\;t\in\Sigma(w)$ & $(\phi\To\psi)\to\ibox(\phi\to\psi)$\\
Finite $\cup$-closure & $s,t\in\Sigma(w)\;\Rightarrow\;s\cup t\in\Sigma(w)$ & $(\aa\To\phi\lori\psi)\to(\aa\To\phi)\lor(\aa\To\psi)$\\
Reflexivity & $wR_\Sigma w$ (i.e., $w\in\bigcup\Sigma(w)$) & $\ibox\aa\to\aa$\\
Non-triviality & $\Sigma(w)\neq\emptyset$ & $\neg{\ibox\bot}$\\
Decreasingness & $wR_\Sigma v\;\Rightarrow\;\Sigma(w)\supseteq\Sigma(v)$ & $(\phi\To\psi)\to\ibox(\phi\To\psi)$\\
Increasingness & $wR_\Sigma v\;\Rightarrow\;\Sigma(w)\subseteq\Sigma(v)$ & $\neg(\phi\To\psi)\to\ibox\neg(\phi\To\psi)$
\end{tabular}
\end{center}
\caption{\label{tab:frame conditions}Some conditions on neighborhood frames and corresponding canonical axioms. In the definition of the conditions, variables are implicitly understood to be universally quantified and $wR_\Sigma v$ is defined as $v\in\bigcup\Sigma(w)$. In the axioms, $\phi$ and $\psi$ stand for arbitrary formulas, while $\aa$ stands for an arbitrary declarative.}
\end{figure}

In this section, we extend the results of the previous section by considering the logic of specific frame classes. 
Six salient conditions on neighborhood frames are listed in Figure~\ref{tab:frame conditions} along with corresponding axiom schemata which, as we now show, are canonical for these conditions, leading to strongly complete axiomatizations.

\begin{theorem} For any subset of the frame conditions in Figure \ref{tab:frame conditions}, a strongly complete axiomatization of \inqnl\ over the resulting class of frames is obtained by adding all instances of the corresponding axiom schemata to the system in the previous section.
\end{theorem}

\smallskip\noindent\emph{Proof.}
Given a subset $\Lambda$ of the frame conditions in Figure \ref{tab:frame conditions}, 
let $\models_\Lambda$ be the consequence relation restricted to the class of in-models whose underlying frames satisfy these conditions, and let $\vdash_\Lambda$ be the relation of derivability obtained by adding to the axiomatization of the previous section all instances of the modal axioms for the conditions in $\Lambda$, as given by Figure \ref{tab:frame conditions}. To show that  $\vdash_\Lambda$ is sound for $\models_\Lambda$ we just need to show that each axiom is valid over models that satisfy the corresponding condition, which is straightforward.

For the converse, we consider a canonical model $M_\Lambda^c$, which is defined as in the previous section, but using $\vdash_\Lambda$ instead of $\vdash$; in other words, worlds in the canonical model are now complete theories of declaratives closed under $\vdash_\Lambda$ (call them $\Lambda$-CTDs for short). Proceeding exactly as in the previous section we can show that if $\Phi\not\vdash_\Lambda\psi$, there is an information state $S$ in $M_\Lambda^c$ that supports all formulas in $\Phi$ but not $\psi$. To obtain completeness, it remains to be shown that the frame of the canonical model $M_\Lambda^c$ satisfies the frame conditions in $\Lambda$, i.e., that each axiom in Figure \ref{tab:frame conditions} is \emph{canonical} for the corresponding condition. Let us consider each axiom in turn.

\medskip\noindent
\textbf{Downward-monotonicity.} Suppose $\vdash_\Lambda$ includes all instances of $(\phi\To\psi)\to\ibox(\phi\to\psi)$. Let $S\in\Sigma^c(\Gamma)$, where $\Gamma$ is a $\Lambda$-CTD, and suppose $T\subseteq S$, $T\neq\emptyset$. Our task is to prove $T\in\Sigma^c(\Gamma)$.
For this, we need to show that whenever $(\phi\To\psi)\in\Gamma$ and $\bigcap T\vdash_\Lambda\phi$, we have $\bigcap T\vdash_\Lambda\psi$. So, suppose $(\phi\To\psi)\in\Gamma$ and $\bigcap T\vdash_\Lambda\phi$. Since $(\phi\To\psi)\to\ibox(\phi\to\psi)$ is an axiom of $\vdash_\Lambda$, $\Gamma$ must contain $\ibox(\phi\to\psi)$, i.e., $(\top\To(\phi\to\psi))$. Since $S\in\Sigma^c(\Gamma)$ and $\bigcap S\vdash_\Lambda\top$, we have $\bigcap S\vdash_\Lambda(\phi\to\psi)$. Since $T\subseteq S$ we have $\bigcap S\subseteq\bigcap T$, and so also $\bigcap T\vdash_\Lambda(\phi\to\psi)$. Since by assumption $\bigcap T\vdash_\Lambda\phi$, we have $\bigcap T\vdash_\Lambda\psi$, as desired.

\medskip\noindent
\textbf{Finite union closure.} Suppose $\vdash_\Lambda$ includes all instances of the schema $(\aa\To\phi\lori\psi)\to(\aa\To\phi)\lor(\aa\To\psi)$, where $\aa$ is a declarative. Let $S,T\in\Sigma^c(\Gamma)$, where $\Gamma$ is a $\Lambda$-CTD. We must prove $S\cup T\in\Sigma^c(\Gamma)$.

To show this, suppose $(\phi\To\psi)\in\Gamma$ and $\bigcap(S\cup T)\vdash_\Lambda\phi$. We need to prove $\bigcap(S\cup T)\vdash_\Lambda\psi$. 
Since $\bigcap(S\cup T)\vdash\phi$, by Lemma \ref{lemma:provable split} we have $\bigcap(S\cup T)\vdash_\Lambda\aa$ for some $\aa\in\R(\phi)$. By Lemma \ref{lemma:resolutions prove formula}, $\aa\vdash_\Lambda\phi$, and so by the (\texttt{CN}) rule, $\vdash_\Lambda(\aa\To\phi)$. Thus, we must have $(\aa\To\phi)\in\Gamma$, and since $(\phi\To\psi)\in\Gamma$, by (\texttt{Trans}) we have $(\aa\To\psi)\in\Gamma$. By Lemma \ref{lemma:provable NF}, this implies $(\aa\To\bb_1\lori\dots\lori\bb_n)\in\Gamma$ where $\{\bb_1,\dots,\bb_n\}=\R(\psi)$. By the schema $(\aa\To\phi\lori\psi)\to(\aa\To\phi)\lor(\aa\To\psi)$, it follows that $(\aa\To\bb_1)\lor\dots\lor(\aa\To\bb_n)\in\Gamma$, and since $\Gamma$ is a complete theory, $(\aa\To\bb)\in\Gamma$ for some particular resolution $\bb\in\R(\psi)$. 

Now, since $S\subseteq S\cup T$ we have $\bigcap(S\cup T)\subseteq \bigcap S$, therefore from $\bigcap(S\cup T)\vdash_\Lambda\aa$ it follows $\bigcap S\vdash_\Lambda\aa$. Since $S\in\Sigma^c(\Gamma)$ and $(\aa\To\bb_i)\in\Gamma$, we must have $\bigcap S\vdash_\Lambda\bb$. Since $\bb$ is a declarative, by Remark \ref{remark:bigcapS} this implies $\bb\in\bigcap S$. With an analogous argument, we can show that $\bb\in\bigcap T$, and therefore $\bb\in\bigcap(S\cup T)$. By Lemma \ref{lemma:resolutions prove formula}, it follows that  $\bigcap(S\cup T)\vdash_\Lambda\psi$, as desired.

\medskip\noindent
\textbf{Reflexivity.} Suppose $\vdash_\Lambda$ includes all instances of the schema $\ibox\aa\to\aa$, where $\aa$ is a declarative. This means that for every $\Lambda$-CTD $\Gamma$, $\ibox\aa\in\Gamma$ implies $\aa\in\Gamma$, which by Lemma \ref{lemma:accessibility} means that $\Gamma R_{\Sigma^c}\Gamma$.

\medskip\noindent
\textbf{Non-triviality.} Suppose $\vdash_\Lambda$ includes the axiom $\neg{\ibox\bot}$ and take any $\Lambda$-CTD $\Gamma$. We must show $\Sigma^c(\Gamma)\neq\emptyset$. We have $\neg{\ibox\bot}\in\Gamma$, that is, $\neg(\top\To\bot)\in\Gamma$ and so $(\top\To\bot)\not\in\Gamma$. The Existence Lemma (Lemma \ref{lemma:existence}) then guarantees the existence of a state $S\in\Sigma^c(\Gamma)$.

\medskip\noindent
\textbf{Decreasingness.} Suppose $\vdash_\Lambda$ includes all instances of the schema $(\phi\To\psi)\to\ibox(\phi\To\psi)$, and suppose $\Gamma R_{\Sigma^c}\Gamma'$ and $S\in\Sigma^c(\Gamma')$. We must show that $S\in\Sigma^c(\Gamma)$. 

For this, suppose $(\phi\To\psi)\in\Gamma$. Then also $\ibox(\phi\To\psi)\in\Gamma$, and therefore by Lemma \ref{lemma:accessibility}, $(\phi\To\psi)\in\Gamma'$. Since $S\in\Sigma(\Gamma')$, we have $\bigcap S\vdash\phi$ implies $\bigcap S\vdash\psi$. Since this holds for every formula $(\phi\To\psi)\in\Gamma$, it follows that $S\in\Sigma^c(\Gamma)$.

\medskip\noindent
\textbf{Increasingness.} Suppose $\vdash_\Lambda$ includes all instances of the schema $\neg(\phi\To\psi)\to\ibox\neg(\phi\To\psi)$, and suppose $\Gamma R_{\Sigma^c}\Gamma'$ and $S\in\Sigma^c(\Gamma)$. We must show that $S\in\Sigma^c(\Gamma')$. 

For this, suppose $(\phi\To\psi)\in\Gamma'$. Then $\neg(\phi\To\psi)\not\in\Gamma'$, and since $\Gamma R_{\Sigma^c}\Gamma'$, by Lemma \ref{lemma:accessibility} we have $\ibox\neg(\phi\To\psi)\not\in\Gamma$. Then from the validity of $\neg(\phi\To\psi)\to\ibox\neg(\phi\To\psi)$, also $\neg(\phi\To\psi)\not\in\Gamma$, and since $\Gamma$ is complete, $(\phi\To\psi)\in\Gamma$. Since $S\in\Sigma^c(\Gamma)$, $\bigcap S\vdash\phi$ implies $\bigcap S\vdash\psi$. Since this holds for every formula $(\phi\To\psi)\in\Gamma'$, it follows that $S\in\Sigma^c(\Gamma')$.

\bigskip
\noindent
The previous theorem should be seen as just the beginning of a systematic study of the inquisitive modal logic of interesting frame classes. In particular, the following table lists some natural frame conditions for which the problem of finding an axiomatization is left for future work.

\begin{center}
\begin{tabular}{lll}
Condition & Definition\\
\hline
$\uparrow$-monotonicity & $s\in \Sigma(w),s\subseteq t\;\Rightarrow\;t\in\Sigma(w)$\\
Convexity & $s,s'\in \Sigma(w),s\subseteq t\subseteq s'\;\Rightarrow\;t\in\Sigma(w)$\\
Finite $\cap$-closure & $s,t\in\Sigma(w)\;\Rightarrow\;s\cap t\in\Sigma(w)$\\
Full $\cup$-closure & $S\subseteq\Sigma(w)\;\Rightarrow\;\bigcup S\in\Sigma(w)$\\
Full $\cap$-closure & $S\subseteq\Sigma(w)\;\Rightarrow\;\bigcap S\in\Sigma(w)$\\
Nestedness & $s,t\in\Sigma(w)\;\Rightarrow\;s\subseteq t\text{ or }t\subseteq s$\\
Weak Centering & $w\in\bigcap\Sigma(w)$
\end{tabular}
\end{center}

\section{Allowing empty neighborhoods}
\label{sec:empty}

We have so far been assuming that neighborhoods are nonempty. What happens if we drop this restriction, allowing the empty set to count as a neighborhood? Our logic $\inqnl$ will not be able discriminate whether or not $\emptyset\in\Sigma(w)$. This is because, by the Empty State property, $\emptyset$ supports every formula, and therefore, the presence or absence of $\emptyset$ in $\Sigma(w)$ can never affect the truth of a modal formula $(\phi\To\psi)$. However, this ``blind spot'' can be removed with a minimal extension of our logic: it suffices to add to our language a modal constant $\odot$, which will count as a modal formula (and thus, as a declarative), with modal depth 1. Its semantics will be truth-conditional, with the following truth conditions:
$$M,w\models\odot\iff \emptyset\in\Sigma(w)$$
Note that when empty neighborhoods are allowed, the formula $\idia\phi$ expresses the existence of a \emph{non-empty} neighborhood supporting $\phi$. The mere existence of a (possibly empty) neighborhood supporting $\phi$ can now be expressed by the formula 
$$\idia_{_\odot}\phi:=\idia\phi\lor\odot$$ 

All the results about the expressivity of $\inqml_{\To}$ proved in Section \ref{sec:expressive power} with respect to inhabited neighborhood models generalize smoothly to the extension $\inqml_{\To\odot}$ with respect to arbitrary neighborhood models. The only change we need to make is to add, in the definition of the characteristic formulas $\chi_{M,w}^n$ for $n\ge 1$, an extra conjunct which is either $\odot$ or $\neg\odot$ depending on whether $\emptyset\in\Sigma(w)$ or not.

The axiomatization given in Section \ref{sec:axiomatization} extends straightforwardly as well. There will be no specific modal axioms ruling $\odot$. However, since $\odot$ counts as a declarative, the propositional axioms ($\texttt{DDN}$) and ($\texttt{Split}$) apply to it, so we have as axioms $\neg\neg\odot\to\odot$ and every instance of $(\odot\to\phi\lori\psi)\to(\odot\to\phi)\lori(\odot\to\psi)$. The definition of the canonical model construction can be extended by letting $\Sigma^c(\Gamma)$ contain $\emptyset$ just in case $\odot\in\Gamma$, which immediately yields the Support Lemma for $\odot$. Everything else remains unchanged.

\section{Relations with instantial neighborhood logic}
\label{sec:inl}

\noindent As discussed in the introduction, instantial neighborhood logic, \inl\ \citep{Benthem:17}, is a modal logic interpreted on neighborhood models which is invariant under the notion of bisimilarity discussed in Section \ref{sec:expressive power}. 
\inl\ is based on a modal language with primitive connectives $\neg$ and $\land$, and where modal formulas have the form $\Box(\rho_1,\dots,\rho_n;\sigma)$ with $n\ge 0$. The semantics is given by a standard definition of truth at a world, where the modal clause is:

\begin{center}
\begin{tabular}{rcl}
$M,w\models_\inl\Box(\rho_1,\dots,\rho_n;\sigma)$ &$\iff$ &$\exists s\in\Sigma(w): (\forall v\in s: M,v\models_\inl\sigma)\text{ and }$\\
&&$\phantom{\exists s\in\Sigma(w):\;}(\forall i\le n\exists v\in s: M,v\models_\inl\rho_i)$
\end{tabular}
\end{center}

\noindent
We will show that \inl\ has the same expressive power as the declarative fragment of $\inqml_{\To\odot}$. We prove this by defining two translations that preserve truth conditions. We define a translation $(\cdot)^*:\L_{\inl}\to\L_{\To\odot}^{!}$ as follows:
\begin{itemize}
\item $p^*=p$
\item $(\neg\sigma)^*=\neg\sigma^*$
\item $(\rho\land\sigma)^*=\rho^*\land\sigma^*$
\item $\Box(\,;\sigma)^*=\idia_{_\odot}\sigma^*$ 
\item $\Box(\rho_1,\dots,\rho_n;\sigma)^*=\neg(\sigma^*\To(\neg\rho_1^*\lori\dots\lori\neg\rho_n^*))$, if $n\ge 1$
\end{itemize}
It is straightforward to verify that this map preserves truth conditions.

\begin{proposition} For every neighborhood model $M$, world $w$, and formula $\sigma\in\L_\inl$, $M,w\models_\inl\sigma\iff M,w\models\sigma^*$.
\end{proposition}

\noindent Translating declaratives of $\inqml_{\To\odot}$ to \inl\ is more tricky. Consider a modal formula $(\phi\To\psi)$: in general, $\phi$ and $\psi$ are not declaratives, so the translation will not be defined on them. Instead, we first compute the resolutions of $\phi$ and $\psi$, which \emph{are} declaratives, and then assemble a translation from the translations of these resolutions. To make this precise, we need a non-standard notion of complexity. Given $\aa,\bb\in\L_!$, we let $\aa\prec\bb$ in case either $\aa$ has lower modal depth than $\bb$, or $\aa$ and $\bb$ have the same modal depth and $\aa$ is a subformula of $\bb$. Clearly, $\prec$ is well-founded and thus suitable for  induction. 
 Now we define a translation $(\cdot)^\star:\L_{\To\odot}^{!}\to\L_{\inl}$ recursively on $\prec$ as follows: 
 \begin{itemize}
\item $p^\star=p$ 
\item $\bot^\star=(p\land\neg p)$ for an arbitrary $p\in\P$
\item $(\aa\land\bb)^\star=\aa^\star\land\bb^\star$
\item $(\aa\to\bb)^\star=\neg(\aa^\star\land\neg\bb^\star)$
\item $\odot^\star=\Box(\,;\bot)$
\item $(\phi\To\psi)^\star=\bigwedge_{i=1}^n\neg\Box(\neg\bb_1^\star,\dots,\neg\bb_m^\star;\aa_i^\star)$\\
where $\{\aa_1,\dots,\aa_n\}=\R(\phi)$ and $\{\bb_1,\dots,\bb_m\}=\R(\psi)$
\end{itemize}
In the last clause, $\aa_i^\star$ is defined since $\aa_i\prec(\phi\To\psi)$: this is because $\aa_i$ has the same modal depth as $\phi$, which is lower that the modal depth of $(\phi\To\psi)$. Similarly, $\bb_j^\star$ is defined since $\bb_j$ has lower modal depth than $(\phi\To\psi)$.

\begin{proposition} Let $\aa\in\L_{\To\odot}^!$ be any declarative in $\inqml_{\To\odot}$. For any neighborhood model $M$ and world $w$ we have $M,w\models\aa\iff M,w\models_\inl\aa^\star$.
\end{proposition}

\begin{proof} By induction on $\prec$. The only interesting step is the case for $\aa=(\phi\To\psi)$. Let $\R(\phi)=\{\aa_1,\dots,\aa_n\}$ and $\R(\psi)=\{\bb_1,\dots,\bb_m\}$. Take any model $M$ and any world $w$ in $M$. We will show that $M,w\not\models(\phi\To\psi)$ iff $M,w\not\models_\inl(\phi\To\psi)^\star$. The second step uses Proposition \ref{prop:normal form}, the fourth uses the fact that declaratives are truth-conditional, and the sixth the induction hypothesis. 
\begin{equation*}
\begin{array}{ll}
&w\not\models(\phi\To\psi)\\
\iff &
\exists s\in\Sigma(w): s\models\phi\text{ and }s\not\models\psi\\
\iff &
\exists s\in\Sigma(w): s\models\Lori_{i\le n}\aa_i\text{ and }s\not\models\Lori_{j\le m}\bb_j\\
\iff &
\exists s\in\Sigma(w): \exists i\le n(s\models\aa_i)\text{ and }\forall j\le m(s\not\models\bb_j)\\
\iff &
\exists s\in\Sigma(w): \exists i\le n(\forall v\in s: v\models\aa_i)\text{ and }\forall j\le m(\exists v\in s: v\not\models\bb_j)\\
\iff &
\exists i\le n\exists s\in\Sigma(w):(\forall v\in s: v\models\aa_i)\text{ and }\forall j\le m(\exists v\in s: v\models\neg\bb_j)\\
\iff &
\exists i\le n\exists s\in\Sigma(w):(\forall v\in s: v\models_\inl\aa_i^\star)\text{ and }\forall j\le m(\exists v\in s: v\models_\inl\neg\bb_j^\star)\\
\iff &
\exists i\le n:w\models_\inl\Box(\neg\bb_{1}^\star,\dots,\neg\bb_m^\star;\aa_i^\star)\\
\iff &
w\not\models_\inl\bigwedge_{i=1}^n\neg\Box(\neg\bb_{1}^\star,\dots,\neg\bb_m^\star;\aa_i^\star)\\
\iff &
w\not\models_\inl(\phi\To\psi)^\star
\end{array}
\end{equation*}
Again, note that we can use the induction hypothesis on $\aa_i$ since $\aa_i\prec(\phi\To\psi)$, and similarly for $\bb_j$. 
\end{proof}

\noindent
A couple of remarks on this translation. First, note that given a formula $\sigma\in\L_\inl$, the size of $\sigma^*$ grows linearly on the size of $\sigma$. By contrast, since the number of resolutions of a formula $\phi\in\L$ grows exponentially in the length of $\phi$ due to the clause for implication, the size of the translation $\aa^\star$ of a formula $\aa\in\L_{\To\odot}^!$ may grow exponentially relative to the size of $\aa$. It seems natural to conjecture that this is inevitable for such a translation, and thus that $\inqml_{\To\odot}$ is exponentially more succinct than \inl. However, we will not try to provide a proof of this conjecture here.
It is also worth noting that this translation strategy, that relies crucially on resolutions, would not be viable in the setting of inquisitive predicate logic, where no analogue of resolutions is available. It is natural to conjecture that a first-order version of $\inqml_{\To\odot}$ would be strictly more expressive than a first-order version of \inl. For instance, one challenge would be to translate the inquisitive modal formula $\ibox(\forall x?Px\to\forall x?Qx)$, which says that in every neighborhood, the extension of $Q$ is functionally determined by the extension of $P$.

\section{Further work}
\label{sec:conclusion}

The work described in this paper can be taken further in various directions. 

First, the mathematical theory of \inqnl\ should be developed further. One immediate goal is to establish the finite model property of \inqnl, and thereby its decidability. Plausibly, this can be achieved by a simple generalization of the filtration method developed by \cite{MaricPerkov:24} for $\inqml_\ibox$ over downward-monotone models. 
It would also be interesting to compare the expressive power of \inqnl\ with that of first-order predicate logic over neighborhood frames, when the latter are encoded as two-sorted structures with distinct domains for worlds and neighborhoods: do we have an analogue of van Benthem's theorem, stating that every bisimulation-invariant first-order formula is equivalent to a modal formula? Again, in previous work \citep{CiardelliOtto:18} this question was answered positively for $\inqml_\ibox$ over downward-monotone models, and it seems that the proof can be generalized smoothly to \inqnl\ over arbitrary models. Alternatively, one could get this result from the analogous result established by \cite{DeGroot:22} for instantial neighborhood logic, via the translations established in Section \ref{sec:inl}.

As mentioned at the end of Section \ref{sec:frame conditions}, the modal logic of many interesting frame classes remains to be worked out. Relatedly, it seems interesting to consider modal formulas at the level of frames, developing a modal correspondence theory for \inqnl. It would be interesting, for instance, to determine what classes of neighborhood frames are definable by modal formulas (in the spirit of the Goldblatt-Thomason theorem for Kripke frames), and to what extent from a modal formula we can cook up a corresponding (two-sorted) first-order formula expressing the same frame condition (in the spirit of Sahlqvist theory).

In a different direction, it would be interesting to explore concrete interpretations of \inqnl, and extensions motivated by such interpretation. For instance, it seems interesting to study a multi-agent version of \inqnl\ interpreted over \emph{concurrent game models} \citep[see][for an overview]{Bulling:16}. In a concurrent game model, worlds represent stages of a dynamic process unfolding over time; at each world $w$, each agent may choose between multiple actions available to them, and the choices of all agents jointly determine the world that will be realized next. A concurrent game model is naturally associated with a multi-agent neighborhood model in the following way: for each world $w$ and each action $x$ available to an agent $a$ at $w$, we may collect the set $\textsf{out}_a(x,w)$ of outcomes that may be realized if $a$ performs $x$ at $w$. We then define a neighborhood structure where $\Sigma_a(w)$ contains as neighborhoods the outcome sets $\textsf{out}_a(x,w)$ for each action $x$ available to $a$ at $w$. In this setting, many formulas of \inqnl\ have an interesting interpretation. Here are some examples:
\begin{itemize}
\item $\idia_a p$ expresses the fact that agent $a$ has an action that guarantees that $p$ will be true at the next stage; this retrieves (for single agents) the semantics of formulas of coalition logic \citep{Pauly:02};
\item $\ibox_a{?p}$ expresses the fact that the action performed by agent $a$ fully determines whether or not $p$ will hold at the next stage---an interesting strategic fact which is provably not expressible in coalition logic;
\item $\neg (p\To_a\neg q)$ expresses the fact that agent $a$ has an action that guarantees the truth of $p$ without precluding the truth of $q$; this sort of property, relevant for the possibility of cooperation between agents with different goals, is central in the \emph{socially friendly coalition logic} developed by \cite{GorankoEnqvist:18}.
\end{itemize}
The multi-agent version of \inqnl\ that we just sketched may be interesting to explore not only for its own sake, but also because it would provide a basis that can be further enriched with proper temporal operators, leading to inquisitive extensions of temporal logics like ATL \citep{Alur:02} and STIT \citep{Belnap:01}.

\bibliographystyle{natbib}
\bibliography{inquisitive}

\end{document}